%% file: 2005-31.tex
\documentclass{gtart_h}  

\input gtoutput

\lognumber{417}
\received{3 February 2004}
\volumenumber{9}\papernumber{31}\volumeyear{2005}
\pagenumbers{1337}{1380}   
\revised{25 July 2005}
\published{30 July 2005}
\accepted{6 May 2004}
\proposed{Haynes Miller}
\seconded{Thomas Goodwillie, Ralph Cohen}

\usepackage{amsmath,amssymb,amscd}

\let\hat\widehat
\let\tilde\widetilde
\let\bar\overline

\numberwithin{equation}{section}

\swapnumbers

\theoremstyle{plain} 
\newtheorem{thm}[equation]{Theorem}

\newtheorem{lem}[equation]{Lemma}
\newtheorem{prop}[equation]{Proposition}

\theoremstyle{definition}
\newtheorem{defn}[equation]{Definition}

\theoremstyle{remark}
\newtheorem{rem}[equation]{Remark}
\newtheorem{aside}[equation]{Aside}

\newtheorem*{VariableNoNum}{{\VariableText}}
\newtheorem{Variable}[equation]{{\VariableText}}

\theoremstyle{definition}
\newtheorem*{VariableNoNumBold}{{\VariableText}}
\newtheorem{VariableBold}[equation]{{\VariableText}}

\newenvironment{numbered}[1]
     {\def\VariableText{{#1}}\begin{Variable}}
     {\end{Variable}}

\newenvironment{NumberedSubSection}[1]
     {\def\VariableText{{\textrm{\textbf{#1}}}}\begin{VariableBold}}
     {\end{VariableBold}}

\def\Math#1{\def\MathString{#1}\futurelet\MathDelim\MathChoose}
\def\MathChoose{\ifmmode\let\MathDo\MathString%
              \else\let\MathDo\MathSkip\fi%
              \MathDo}
\def\MathSkip{\ifx\MathDelim/\def\MathDo{$\MathString$\EatOne}%
              \else\def\MathDo{$\MathString$}\fi%
              \MathDo}

\def\Text#1{\def\TextString{#1}\futurelet\TextDelim\TextSkip}
\def\TextSkip{\ifx\TextDelim/\def\TextDo{\TextString\EatOne}%
              \else\let\TextDo\TextString\fi%
              \TextDo}
\def\EatOne#1{}

\def\SkipToEndScan#1\EndScan{}
\def\Scan#1#2#3{\ifx#1#2#3\expandafter\SkipToEndScan\fi\Scan#1}
\def\Upper#1{%
\Scan#1aAbBcCdDeEfFgGhHiIjJkKlLmMnNoOpPqQrRsStTuUvVwWxXyYzZ#1#1\EndScan}
\def\Phrase#1 #2/#3/#4=#5 #6/#7/#8.{%
\expandafter\edef\csname#2#3\endcsname{\noexpand\Text{#6#7}}
\expandafter\edef\csname\Upper#2#3\endcsname{\noexpand\Text{\Upper#6#7}}
\expandafter\edef\csname#1#2#3\endcsname{\noexpand\Text{#5 #6#7}}
\expandafter\edef\csname\Upper#1#2#3\endcsname{\noexpand\Text{\Upper#5 #6#7}}
\expandafter\edef\csname#2#4\endcsname{\noexpand\Text{#6#8}}
\expandafter\edef\csname\Upper#2#4\endcsname{\noexpand\Text{\Upper#6#8}}
}
\newcommand{\R}{\Math{\mathbb{R}}}
\newcommand{\C}{\mathbb{C}}
\newcommand{\Z}{\Math{\mathbb{Z}}}
\newcommand{\Q}{\Math{\mathbb{Q}}}
\newcommand{\F}{\mathbb{F}}

\newcommand{\Tensor}{\otimes}

\newcommand{\RightArrow}[1]{\xrightarrow{#1}} 

\newcommand{\GL}{\operatorname{GL}}

\newcommand{\HH}{\operatorname{H}}

\newcommand{\Hom}{\operatorname{Hom}}
\newcommand{\Aut}{\operatorname{Aut}}
\newcommand{\Map}{\operatorname{Map}}

\newcommand{\hhh}{\operatorname{h}\!}
\newcommand{\thhh}{\tilde{\operatorname{h}}}
\def\HomotopyOrbit#1on#2/{\ensuremath{#2_{\hhh#1}}}
\def\RedHomotopyOrbit#1on#2/{\ensuremath{#2_{\thhh#1}}}

\newcommand{\iso}{\cong}
\newcommand{\diag}{\operatorname{diag}}

\Phrase a ts//s=a marking//s.
\Phrase a marking//s=a marking//s.
\Phrase a markedtor/us/i=a marked tor/us/i.
\Phrase a markedreflectiontor/us/i=a marked reflection tor/us/i.
\Phrase a markedreflectionlattice//s=a marked reflection lattice//s.
\Phrase a markedreflectionstructure//s=a marked reflection structure//s.

\newcommand{\Ztinfty}{\Math{\Z/2^\infty}}
\newcommand{\weq}{\sim}
\newcommand{\Zt}{\Math{\Z_2}}
\newcommand{\Qt}{\Math{\Q_2}}
\newcommand{\Exten}{\nu}

\newcommand{\Ztwo}{\Math{\Z_2}}
\newcommand{\Ftwo}{\Math{\F_2}}
\Phrase a ttor/us/i=a discrete tor/us/i.
\Phrase a tlattice//s=a complete lattice//s.
\Phrase an rtlattice//s=a complete reflection lattice//s.
\Phrase a mrtlattice//s=a marked complete reflection lattice//s.
\Phrase an rttor/us/i=a discrete reflection tor/us/i.
\Phrase a mrttor/us/i=a marked discrete reflection tor/us/i.
\newcommand{\td}{\Text{$2$--discrete}}

\newcommand{\DiscreteMarker}[1]{\breve{#1}}
\newcommand{\CompleteMarker}[1]{\hat{#1}}
\newcommand{\Td}{\Math{\DiscreteMarker T}}
\newcommand{\Sd}{\Math{\DiscreteMarker S}}
\newcommand{\Scl}{\Math{\bar S}}
\newcommand{\Nd}{\Math{\DiscreteMarker N}}
\newcommand{\Lc}{\Math{\CompleteMarker L}}

\newcommand{\Tc}{\Math{\CompleteMarker T}}
\newcommand{\LD}{\Math{L_{\Delta}}}
\newcommand{\WD}{\Math{W_{\Delta}}}
\newcommand{\cTD}{\Math{\mathcal{T}_{\Delta}}}
\newcommand{\Qtwo}{\Math{\Q_2}}
\newcommand{\RG}[1]{G(#1)}

\newcommand{\RX}[1]{X(#1)}
\newcommand{\RXA}[1]{X_A(#1)}

\newcommand{\SU}{\Math{\operatorname{SU}}}
\newcommand{\SO}{\Math{\operatorname{SO}}}
\newcommand{\U}{\Math{\operatorname{U}}}
\newcommand{\SUc}{\Math{\hat{\operatorname{S}}\!\operatorname{U}(2)}}
\newcommand{\SOc}{\Math{\hat{\operatorname{S}}\!\operatorname{O}(3)}}
\newcommand{\Uc}{\Math{\hat{\operatorname{U}}(2)}}

\newcommand{\prodd}{\operatorname{prod}}
\newcommand{\Torus}{\operatorname{T}}
\newcommand{\dual}{^\#}
\newcommand{\im}{\operatorname{im}}
\newcommand{\RE}{\rho}
\newcommand{\mcT}{\mathcal{T}}
\newcommand{\Tits}{\tau}

\newcommand{\Iset}{\Math{I}}
\newcommand{\Igroup}{\Math{\mathbf{I}}}
\newcommand{\Iel}{\Math{\mathbf{i}}}
\newcommand{\Jel}{\Math{\mathbf{j}}}
\newcommand{\Norm}{\nu}
\newcommand{\Normd}{\breve\nu}
\newcommand{\refl}{\sigma}
\newcommand{\crefl}{\tau} 
\newcommand{\refls}{\Sigma}
\newcommand{\srefl}{s}
\newcommand{\srefls}{S}
\newcommand{\dTorus}{\breve{\operatorname{T}}}
\newcommand{\cTorus}{\hat{\operatorname{T}}}
\newcommand{\pperp}{^\perp}

\newcommand{\DI}{\operatorname{DI}}
\begin{document}

\title{Normalizers of tori}
\author{W\,G Dwyer\\C\,W Wilkerson}
\asciiauthors{WG Dwyer and CW Wilkerson}
\coverauthors{W\noexpand\thinspace G Dwyer\\C\noexpand\thinspace W Wilkerson}

\address{Department of Mathematics, University of Notre Dame\\Notre Dame,
Indiana 46556, USA}
\address{Department of Mathematics, Purdue University\\West Lafayette,
Indiana 47907, USA}
\asciiaddress{Department of Mathematics, University of Notre Dame\\Notre Dame,
Indiana 46556, USA\\and\\Department of Mathematics, Purdue University\\West Lafayette,
Indiana 47907, USA}

\gtemail{\mailto{dwyer.1@nd.edu}, \mailto{wilker@math.purdue.edu}}
\asciiemail{dwyer.1@nd.edu, wilker@math.purdue.edu}


\begin{abstract}
We determine the groups which can appear as the normalizer of a
maximal torus in a connected
$2$--compact group. The technique depends on using ideas of Tits to give
a novel description of the normalizer of the torus in a connected
compact Lie group, and then showing that this description can  be
extended to the $2$--compact case.
\end{abstract}
\asciiabstract{%
We determine the groups which can appear as the normalizer of a
maximal torus in a connected 2-compact group. The technique depends
on using ideas of Tits to give a novel description of the normalizer
of the torus in a connected compact Lie group, and then showing that
this description can be extended to the 2-compact case.}

\maketitlepage

\section{Introduction}

Suppose that $G$ is a connected compact Lie group and that $T\subset
G$ is a maximal torus, or in other words a maximal connected
abelian subgroup. The normalizer $NT$ of $T$ lies in a
short exact sequence
\begin{equation}\label{BasicExtension}
  1 \to T \to NT \to W \to 1
\end{equation}
in which $W$ is a finite group called the \emph{Weyl group} of $G$.
In this paper we use ideas of Tits \cite{rTits} to give a particularly
simple description 
of the groups which appear as such an $NT$ (see Proposition~\ref{NTForG}). This leads to
an analogous determination of the groups which appear as the normalizer
$N\Td$ of a maximal $2$--discrete torus $\Td$ in a connected $2$--compact
group (see Section~\ref{Notation}).

In the connected compact Lie group case, $NT$ determines $G$ up to isomorphism;
see \cite{rCurtisEtAl} for  semisimple groups, and \cite{rOsse} 
or \cite{rNotbohmFunctor} in general. In listing the possible $NT$s we are
thus giving an alternative approach to the classification of connected
compact Lie groups. In contrast, it is not known that the
normalizer of a maximal $2$--discrete torus in a connected $2$--compact
group $X$ determines $X$ up to equivalence. However, this seems likely
to be true \cite{rNotbohmOrthogonal,rMollerDetermined,rOddClassification},
and we hope that the results of this paper will eventually contribute
to a classification of connected $2$--compact groups.

\begin{NumberedSubSection}{Compact Lie groups}
A glance at \eqref{BasicExtension} reveals that $NT$ is determined up to
isomorphism by three ingredients:
\begin{enumerate}
\item the torus $T$, \label{TorusIngredient}
\item the  finite group $W$,  with its conjugation action  on
  $T$, and \label{WeylIngredient}
\item an extension class $k\in H^2(W;T)$. \label{ExtensionIngredient}
\end{enumerate}
Let $\Torus(1)$ denote the circle $\R/\Z$.
The torus $T$ is a product $\Torus(1)^r$ of circles, and so is
determined by the number $r$, which is the \emph{rank} of $T$ (or of
$G$). The natural map $\Aut(T)\to\Aut(\pi_1T)\iso\GL(r,\Z)$ is an
isomorphism, and the resulting conjugation homomorphism $W\to
GL(r,\Z)$ embeds $W$ as a finite subgroup of $\GL(r,\Z)$ generated by
reflections.  Here a \emph{reflection} is a matrix which is conjugate
in $\GL(r,\Q)$ to the diagonal matrix $\diag(-1,1,\ldots,1)$. 
This takes care of ingredients \eqref{TorusIngredient} and
\eqref{WeylIngredient}.  Ingredient \eqref{ExtensionIngredient} is
more problematical. The most direct way to approach
\eqref{ExtensionIngredient} would be to list all of the elements of
$H^2(W;T)$ and then point to ones which correspond to $NT$ extensions,
but it would be a daunting task to make all of these cohomology
calculations (one for each finite group acting on a torus by
reflections) and then name the resulting elements in a way which would
make it convenient to single out the ones giving an $NT$.

We take another approach, which is based on \cite{rTits}. Suppose that
$T$ is a torus and $W$ is a finite group of automorphisms of $T$
generated by reflections. \emph{\Ats/} for $(T,W)$ is defined to be a
collection $\{h_\refl\}$ of elements of $T$, one for each reflection
$\refl$ in $W$, which satisfy some simple conditions
(see Definition~\ref{DTorusMarking}). The triple $(T,W,\{h_\refl\})$ is called a
\emph{marked reflection torus}.  The following proposition is less
deep than it may seem; it is derived from a simple algebraic
correspondence between marked reflection tori and classical root
systems (see Section~\ref{CMenagerie}).

\begin{prop}[See Remark~\ref{MarkingsAndClasicalRoots}]\label{ClassifyG}
   Suppose that $G$ is a connected
  compact Lie group with maximal torus $T$ and Weyl group $W$. Then
  $G$ determines a natural marking $\{h_\refl\}$ for the pair $(T,W)$,
  and the assignment $G\mapsto (T,W,\{h_\refl\})$ gives a bijection
  between isomorphism classes of connected compact Lie groups and
  isomorphism classes of marked reflection tori.
\end{prop}

To each such marked reflection torus we associate a group
$\Exten(T,W,\{h_\refl\})$, called the \emph{normalizer extension of
  $W$ by $T$} (see Definition~\ref{DefineNormalizerExt}), which lies in a short exact
sequence
\begin{equation}\label{MarkedTorusExtension}
1 \to T \to \Norm(T,W,\{h_\refl\})\to W \to 1.
\end{equation}
The following result is essentially due to Tits \cite{rTits}.

\begin{prop}[See Theorem~\ref{IdentifyTorusNormalizer}]
  \label{NTFromMarked}
  Suppose that $G$ is a connected compact Lie group
  and that
  $(T,W,\{h_\refl\})$ is the corresponding marked
  reflection torus. Then the normalizer of $T$ in $G$ is
  isomorphic to $\Norm(T,W,\{h_\refl\})$.
\end{prop}

This leads to our characterization of $NT$s.

\begin{prop}\label{NTForG}
  Suppose that $N$ is an extension of a finite group $W$ by a
  torus~$T$.  Then $N$ is isomorphic to the normalizer of a torus in a
  connected compact Lie group if and only if
  \begin{itemize}
  \item the conjugation action of $W$ on $T$ expresses $W$ as a group
    of automorphisms of $T$ generated by reflections, and
  \item there exists a marking $\{h_\refl\}$ of $(T,W)$ such that $N$ is
    isomorphic to $\Norm(T,W,\{h_\refl\})$.
  \end{itemize}
\end{prop}

In fact the marking $\{h_\refl\}$  in Proposition~\ref{NTForG} is unique, if it
exists (see Remark~\ref{MarkingsAndClasicalRoots}).

\end{NumberedSubSection}

\begin{NumberedSubSection}{$2$--compact groups}
  A \emph{$2$--compact group} is by definition a pair $(X, BX)$ of
  spaces together with an equivalence $X\weq\Omega BX$, where
  $H_*(X;\Z/2)$ is finite and $BX$ is $2$--complete. Usually we refer to
  $X$ itself as the $2$--compact group and leave $BX$ understood. Any
  such $X$ has a maximal \emph{$2$--discrete torus} $\Td$
  (see Section~\ref{Notation}), which is an ordinary discrete group of the form
  $(\Ztinfty)^r$ for some $r>0$.  There is also an associated
  normalizer $N\Td$ which lies in a short exact sequence
  \[
   1 \to \Td \to N\Td \to W \to 1
  \]
  in which $W$ is a finite group called the \emph{Weyl group} of~$X$. The automorphism group of $\Td$ is
  isomorphic to $\GL(r,\Zt)$, and again the conjugation action of $W$
  on $\Td$ embeds $W$ as a subgroup of $\Aut(\Td)$ generated by
  reflections. In this context a \emph{reflection} is an element of
  $\GL(r,\Zt)$ which is conjugate in $\GL(r,\Qt)$ to
  $\diag(-1,1,\ldots,1)$. As above, we define a \emph{marking} for
  $(\Td,W)$ to be a suitable (see Definition~\ref{tDefineTorusMarking}) collection of
  elements $\{h_\refl\}$ in $\Td$, and we associate to each such
  marking a group $\Norm(\Td,W,\{h_\refl\})$, called the
  \emph{normalizer extension of $W$ by $\Td$}
  (see Definition~\ref{tDefineNormalizerExt}), which lies in a short exact sequence
  \begin{equation}\label{FundamentalExtension}
   1 \to \Td \to \Norm(\Td,W,\{h_\refl\})\to W\to 1. 
  \end{equation}
  The structure $(\Td,W,\{h_\refl\})$ goes by the name of \emph{marked
    $2$--discrete reflection torus}.  Our main results are the
  following ones.

\begin{prop}[See Lemma~\ref{LXMarkings}]\label{ClassifyX} 
  Suppose that $X$ is a connected $2$--compact group with maximal
  $2$--discrete torus $\Td$ and Weyl group $W$. Then $X$ determines a
  natural marking $\{h_\refl\}$ for $(\Td,W)$, and thus a marked $2$--discrete
  reflection torus $(\Td,W,\{h_\refl\})$.
\end{prop}

\begin{prop}[See Proposition~\ref{IdentifyNT}]\label{NTForX}
  Suppose that $X$ is a connected $2$--compact group
  and that $(\Td,W,\{h_\refl\})$ is the
  corresponding (see Proposition~\ref{ClassifyX}) marked $2$--discrete reflection torus.
  Then the normalizer of $\Td$ in $X$ is isomorphic to
  $\Norm(\Td,W,\{h_\refl\})$. 
\end{prop}

It turns out that there are not very many marked $2$--discrete reflection
tori to work with.  If $(T,W,\{h_\refl\})$ is a marked reflection torus,
there is an associated marked $2$--discrete reflection torus
$(T^\delta, W, \{h_{\refl}\})$, where
\begin{equation}\label{TwoPrimary}
      T^\delta=\{x\in T\mid 2^kx=0\text{ for }k\gg 0\}.
\end{equation}
Say that a $2$--discrete marked reflection
torus is \emph{of Coxeter type} if it is derived in this way from an
ordinary marked reflection torus.
  Let $\cTD$
denote the marked $2$--discrete reflection torus derived from the
exceptional $2$--compact group $\DI(4)$ \cite{rDWnew,rNotbohmDIFour}.  There is  an evident notion of cartesian product for marked
$2$--discrete reflection tori (see Section~\ref{CtLattices}).

\begin{prop}[See Propositions~\ref{tInterpretMarkedTorus} and \ref{ClassifyCompleteLattices}]
\label{DoubleProduct}  
Any marked $2$--discrete reflection torus can be written as a product
  $\mathcal T_1\times \mathcal T_2$, where $\mathcal T_1$ is of
  Coxeter type and $\mathcal T_2$ is a product of copies of $\cTD$.
\end{prop}

It follows from Proposition~\ref{DoubleProduct} that any marked $2$--discrete
reflection torus arises from some $2$--compact group: those of
Coxeter type from $2$--completions of connected compact Lie groups
(Section~\ref{CPassingToTwo} and Lemma~\ref{SameRootSubgroups}), and the
others from products of copies of $\DI(4)$. This can be used to give
a characterization parallel to Proposition~\ref{NTForG} of normalizers of maximal
$2$--discrete tori in connected $2$--compact groups.

\begin{rem}
  Proposition \ref{DoubleProduct} is in line with the conjecture that
  any connected $2$--compact group $X$ can be written as $X_1\times
  X_2$, where $X_1$ is the $2$--completion of a connected compact Lie
  group and $X_2$ is a product of copies of $\DI(4)$ (see
  Proposition~\ref{SplitIntoCases}). Taken together, Propositions~\ref{NTForX} and
  \ref{DoubleProduct} imply that the conjecture is true as far as
  normalizers of maximal $2$--discrete tori are concerned.
\end{rem}

\begin{rem}
  If $p$ is an odd prime, the normalizer of a maximal $p$--discrete
  torus in a connected $p$--compact group is always the semidirect
  product of the torus and the Weyl group \cite{rAndersen}. This
  explains why we concentrate on $p=2$; for odd $p$ the extension we
  study is trivial. The techniques of this paper work to some extent
  for odd primes and can be
  used to show that the extension is trivial when the
  Weyl group is a Coxeter group; most of the remaining cases can be
  handled with the general methods of \cite{rDWCoxeter}.
\end{rem}

\end{NumberedSubSection}

\begin{NumberedSubSection}{Organization of the paper}
  Section~\ref{CMenagerie} describes various ways of interpreting a root
  system, one of which is in terms of a \emph{marked reflection
    torus}.  Section~\ref{CNormalizerExt} shows how a marked
  reflection torus $(T,W,\{h_\refl\})$ gives rise to a
  \emph{normalizer extension} $\Norm(W)$ of $W$ by $T$; the
  normalizer extension is defined with the help of a \emph{reflection
    extension} $\RE(W)$ of $W$.
  In Section~\ref{CExtended} we find generators and
  relations for $\RE(W)$ by identifying this group with an extension
  $\Tits(W)$ of $W$ constructed by Tits. Section~\ref{CCompactLie}
  exploits these generators and relations to show that  the
  normalizer $NT$ of a maximal torus in a connected compact Lie group
  $G$ is isomorphic to the normalizer extension obtained from the marked reflection torus (equivalently root system)
  associated to $G$. The arguments are due to Tits \cite{rTits},
  although we put them in a different context and present them in a 
  way 
  that can be generalized to $2$--compact groups. 
  Section~\ref{CMenagerieTwo} generalizes the
  material from Sections~\ref{CMenagerie} and \ref{CNormalizerExt} to a
  $2$--primary setting, and Section~\ref{CtLattices} goes on to give a
  classification of marked $2$--discrete reflection tori.
  Section~\ref{CPassingToTwo} gathers together some technical information
  about
  $2$--completions of connected compact Lie groups. Finally,
  Section~\ref{CTCGroups} proves the main results about the normalizer of
  a maximal $2$--discrete torus in a connected $2$--compact group;
  the approach
  depends in part on a technical lemma, which is proved in
  Section~\ref{CCentralizer}.
\end{NumberedSubSection}

\begin{NumberedSubSection}{Notation and terminology}\label{Notation}
  By a \emph{maximal $2$--discrete torus} $\Td$ in a $2$--compact group
  $X$, we mean the discrete approximation \cite[6.4]{rDWmethods} to a
  maximal $2$--compact torus $\Tc$ in $X$ \cite[8.9]{rDWmethods}
  \cite[2.15]{rDWcenter}. The \emph{normalizer} of $\Td$ is the discrete
  approximation in the sense of \cite[3.12]{rDWcenter} to the
  normalizer of $\Tc$ \cite[9.8]{rDWmethods} (Aside~\ref{DiscussDiscreteTorus}). We take for granted the
  basic properties of $2$--compact groups (see for instance the survey
  articles \cite{rLannesSurvey,rMollerSurvey,rNotbohmSurvey}), although
  we sometimes recall the definitions
  to help orient the reader. Since we only work with $p$--compact
  groups for $p=2$, we sometimes abbreviate \emph{$2$--complete} to
  \emph{complete} and \emph{$2$--discrete} to \emph{discrete}.

  Given an abelian group $A$ and an involution $\refl$ on $A$, we
  write $A^-(\refl)$ for the kernel of $1+\refl$ on $A$, and
  $A^+(\refl)$ the kernel of $1-\refl$.

The following equations give the relationship between our notions of
\emph{root system} (Definition~\ref{DefineOurRootSystems}) and \emph{geometric
  root system} (see Section~\ref{TwoRootSystems})  and similar
notions of  Tits, Serre, and Bourbaki.
\[
\begin{aligned}
\text{\emph{root system}} &= \text{syst\`eme de racines \cite{rTits}}\\& =
\text{root data \cite{rSpringerLAG}}\\
    &=\text{diagramme radiciel r\'eduit \cite[\S4, Sect.~8]{rBourbakiNine}}\\
\text{\emph{geometric root system}}&= \text{syst\`eme de racines
  \cite{rSerre} \cite{rBourbaki}} 
\end{aligned}
\]
In a nutshell, for us \emph{root systems} are integral objects, while
\emph{geometric root systems} are real objects.
It is a classical result \cite{rBourbaki} that geometric root systems
classify connected compact semisimple Lie groups up to isogeny;
equivalently, they classify simply connected compact Lie groups (or
center--free connected compact Lie groups) up to isomorphism.  A root
system (in our sense) is an integral analog of a geometric root
system, an analog which has the additional feature of allowing for the existence
of central tori; root systems classify arbitrary connected compact Lie
groups up to isomorphism.  The existence of a natural correspondence
between root systems and connected compact Lie groups is documented in
\cite[\S4, Sect.~9]{rBourbakiNine}.  This correspondence can be
derived fairly easily from the classical correspondence between simple
geometric root systems and simply connected simple compact Lie groups. 
However, there is an indirect path to
the correspondence that goes through the intermediate category of
connected reductive linear algebraic groups over $\mathbb C$
(see Section~\ref{TwoRootSystems}).
\end{NumberedSubSection}

The authors would like to thank Kasper Andersen, Jesper Grodal, and
the referee for their comments and suggestions; they are particularly
grateful to T\,A Springer for pointing out some references to the
literature.

\section{Root systems, lattices and tori}\label{CMenagerie}
\label{CRootSystems}

In this section we show that three different structures are
equivalent: \emph{root systems}, \emph{marked reflection lattices},
and \emph{marked reflection tori}.  The conclusion is that a root system
amounts to a reflection lattice $L$ together with a bit of extra
marking data for each reflection; this marking data is particularly
easy to specify in terms of the torus $\Torus(1)\Tensor L$.

\begin{NumberedSubSection}{Root systems}
The following definition is adapted from \cite[4.1]{rTits}.
A \emph{lattice} $L$ (over \Z) is a free abelian group of finite rank.

  \begin{defn}\label{DefineRootSystem}\label{DefineOurRootSystems}
    A \emph{root system} in a lattice $L$ consists of a finite subset
    $R$ of $L$ (the set of \emph{roots}) and for each $ r \in R$
    a homomorphism $n_ r \co  L\to\Z$ (the \emph{coroot} associated to
    $ r $) such that following axioms are satisfied:
    \begin{itemize}
    \item[(R1)] taken together, $R$ and 
      $\bigcap_ r  \ker(n_ r )$ span $\Q\Tensor L$;
    \item[(R2)] for each $ r \in R$, $n_{ r }( r )=-2$; 
    \item[(R3)] if $ r \in R$, $k\in \Z$, and $k r \in R$, then
      $k=\pm1$;
    \item[(R4)] if $ r ,t\in R$, then $t+n_ r (t) r \in R$.
    \end{itemize}
  \end{defn}

  \begin{rem}\label{PointOutDifference}
    This structure is called a \emph{syst\`eme de racines} in
   \cite{rTits} and a \emph{diagramme radiciel r\`eduit} in
   \cite{rBourbakiNine}. It is not the conventional \emph{syst\`eme de
     racines} of Serre \cite[Section~5]{rSerre} or Bourbaki
    \cite[VI]{rBourbaki}, which we will call a \emph{geometric root
      system}.  Basically, a root system in our sense is an integral
    form of a geometric root system; see Section~\ref{TwoRootSystems} below.
  \end{rem}
\end{NumberedSubSection}

\begin{NumberedSubSection}{Marked Reflection Lattices}
  A \emph{reflection} on a lattice $L$ is an automorphism $\refl$
  such that $\refl$ is conjugate in  $\Aut(\Q\Tensor
  L)\cong\GL(r,\Q)$ to a diagonal matrix $\diag(-1,1,\ldots,1)$.

\begin{defn}\label{DefineMarking}
  Suppose that $\refl$ is a reflection on a lattice $L$. A
  \emph{strict marking} for $\refl$ is a pair $(b,\beta)$, where
  $b\in L$ and $\beta\co L\to\Z$ is a homomorphism such that for any
  $x\in L$
  \begin{equation}\label{MarkingEquation}
     \refl(x) = x + \beta(x) b.
  \end{equation}
  Two strict markings $(b,\beta)$ and $(b',\beta')$ are
  \emph{equivalent} if $(b,\beta)=\pm(b',\beta')$.  A \emph{marking}
  for $\refl$ is an equivalence class of strict markings.
\end{defn}

Suppose that $L$ is a lattice, $\refl$ is a reflection on $L$ with
marking $\pm(b,\beta)$, and $w$ is an automorphism of $L$. We let
$w\cdot(b,\beta)$ denote the marking for $w\refl w^{-1}$ given by
$(w(b),\beta\circ w^{-1})$.

\begin{defn}\label{DefineMarkedLattice}
  A \emph{reflection lattice} is a lattice $L$ together with a finite
  subgroup $W$ of $\Aut(L)$ which is generated
  by the reflections it contains. A \emph{marked reflection lattice}
  is a reflection lattice $(L,W)$ together with markings
  $\pm(b_\refl,\beta_\refl)$, one for each reflection $\refl$ in
  $W$, such that for
  $w\in W$,
  \[
      w\cdot (b_\refl, \beta_\refl)= \pm(b_{w\refl w^{-1}},
        \beta_{w\refl w^{-1}}).
  \]
\end{defn}

\end{NumberedSubSection}

\begin{NumberedSubSection}{Marked reflection tori}

  Let $\Torus(n)$ denote the $n$--torus $(\R/\Z)^n$. A \emph{torus} is
  a compact Lie group isomorphic to $\Torus(n)$ for some $n$. Any
  torus $T$ gives a lattice $\pi_1T$; conversely, a lattice $L$ gives
  a torus $\Torus(1)\Tensor L$. These two constructions are inverse to
  one another up to natural isomorphism, and induce an equivalence
  between the category of tori (with continuous homomorphisms) and the
  category of lattices.

\begin{rem}
  We write the group operation in a torus additively.
  If $T$ is a torus, the lattice $\pi_1T$ can also be described as the
  group $\Hom(\Torus(1),T)$ of continuous homomorphisms from the $1$--torus
  to~$T$ (these are \emph{cocharacters} of $T$).
\end{rem}

\begin{defn}
  An automorphism $\refl$ of $T$ is said to be a \emph{reflection}
  (respectively,  \emph{trivial mod~$2$}) if the induced automorphism of $\pi_1T$ is
  a reflection (respectively, trivial mod~$2$).  An element $x$ of $T$ is said
  to be \emph{strongly $\refl$--negative} if $x$ lies in the identity
  component $T^-_0(\refl)$ of $T^-(\refl)$ (see Section~\ref{Notation}), or, in
  other words, if $x$ lies in the identity component of the group
  \[
      T^-(\refl)=\{x\in T\mid \refl(x)=-x\}.
  \] 
   
\end{defn}

\begin{rem}\label{DefineStronglyNegative}
  A reflection $\refl$ on $T$  is trivial mod~$2$ if and
  only if $\refl$ acts as the identity on ${}_2T=\{x\in T\mid
  2x=0\}$. If $T=\Torus(1)\Tensor L$, then $x\in T$ is strongly
  $\refl$--negative if and only if $x\in \Torus(1)\Tensor L^-(\refl)$.
\end{rem}

\begin{defn}\label{DTorusMarking}
  Suppose that $\refl$ is a reflection on a torus $T$. A
  \emph{marking} for $\refl$ is an element $h\in T$ such that
  \begin{enumerate}
  \item $h$ is strongly $\refl$--negative (Remark~\ref{DefineStronglyNegative}),
  \item $2h=0$, and
  \item $h\ne0$ if $\refl$ is nontrivial mod~$2$.
  \end{enumerate}
\end{defn}

\begin{defn}\label{DMarkedReflectionTorus}
  A \emph{reflection torus} is a torus $T$ together with a finite
  subgroup $W$ of $\Aut(T)$ which is generated by the reflections it
  contains. A \emph{marked reflection torus} is a reflection torus
  $(T,W)$ together with markings $h_\refl\in T$, one for each
  reflection $\refl$ in $W$, such that for $w\in W$,
  \[
          h_{w\refl w^{-1}} = w(h_\refl).
  \]
\end{defn}

\end{NumberedSubSection}

\begin{rem}\label{CoxeterGeneration}
  If $(L,W)$ is a reflection lattice, or $(T,W)$ is a reflection
  torus, then $W$ is a classical Weyl group \cite[2.9]{rHump} and
  hence a Coxeter group \cite[Section~1]{rHump}. Let
  $\srefls=\{\srefl_1,\ldots,\srefl_\ell\}$ be a set of simple reflections in $W$
  \cite[1.3]{rHump}, and $m_{i,j}$ the order of $\srefl_i\srefl_j$, so
  that $\{m_{i,j}\}$ is the Coxeter matrix associated to $\srefls$. By
  \cite[1.9]{rHump}, $W$ is generated by $\srefls$ subject only to the
  relations $(\srefl_i\srefl_j)^{m_{i,j}}=1$. 
\end{rem}

\begin{NumberedSubSection}{Equivalence between the three structures}

The goal of this section is to prove the following two statements.

\begin{prop}\label{InterpretRootSystem}
    If $L$ is a lattice, then there is a natural bijection between
    marked reflection structures on $L$ and root systems in $L$.
\end{prop}

\begin{prop}\label{InterpretMarkedTorus}
  If $L$ is a lattice with associated torus $T=\Torus(1)\Tensor L$,
  then there is a natural bijection between marked reflection
  structures on $L$ and marked reflection structures on $T$.
\end{prop}

\begin{rem}\label{HowManyMarkings} It is not hard to enumerate the
  possible markings for a reflection $\refl$ on $L$.  Write
  $\ker(1+\refl)$ and $\im(1-\refl)$ for the kernel and image of the
  indicated endomorphisms of $L$; since $\refl$ is a reflection, both
  of these groups are infinite cyclic.  If $b\in L$ is an element such
  that $\im(1-\refl)\subset\langle b\rangle$, then for every $x\in L$
  there is a unique integer $ \beta(x)$ such that
  $\refl(x)=x+\beta(x)b$; it is easy to check that $\beta\co L\to \Z$ is
  a homomorphism, and so $(b,\beta)$ is a strict marking. In
  particular, strict markings for $\refl$ correspond bijectively
  to such elements~$b$, and markings to such elements $b$ taken up to
  sign. Any such $b$ lies in $\ker(1+\refl)$.  It is easy to see that
  there are inclusions
  \[
    2 \ker(1+\refl) \subset \im(1-\refl)\subset \ker(1+\refl).
  \]
  where the left hand inclusion is an equality if and only if $\refl$
  reduces to the identity mod~$2$.
  Let $b_0$ be a generator of $\ker(1+\refl)$.  Then
  $\refl$ has exactly two markings (determined by $\pm b_0$ and
  $\pm2b_0$) if $\refl$ is trivial mod~$2$, and one marking
  (determined by $\pm b_0$) otherwise.
\end{rem}

\begin{proof}[Proof of Proposition~\ref{InterpretRootSystem}]
  Suppose to begin with that $(L,W)$ is a marked reflection lattice.
  For each reflection $\refl\in W$, let $\pm(b_\refl,\beta_\refl)$
  be the corresponding marking, so that
  \begin{equation}\label{RootMarking}
       \refl(x)=x+\beta_\refl(x)b_\refl.
  \end{equation}
  The set of roots in the corresponding root system is given by
  \[
  R=\{ \pm b_\refl\mid \refl\in W\text{ a reflection}\} .
  \]
  For each root $r\in R$, the coroot $n_r$ is $\beta_\refl$ if
  $r=b_\refl$ and $-\beta_\refl$ if $r=-b_\refl$. It remains to
  check axioms (R1)-(R4).  For (R1), note that $\bigcap_r\ker(n_r)$ is
  the fixed point set $L^W$ of the action of $W$ on $L$, and observe
  that for any $x\in L$ there is an expression
  \[
   (\#W)x = \sum_{w\in W} w(x) + \sum_{w\in W} (x-w(x)),
  \]
  with the first summand on the right in $L^W$ and the second in the
  span of $R$. Condition (R2) comes from combining \eqref{RootMarking}
  with the equality $\refl^2(x)=x$.  For each $w\in W$ and reflection
  $\refl\in W$, $w(b_\refl\vphantom{w^{-1}})=\pm b_{w\refl w^{-1}}$
  (Definition~\ref{DefineMarkedLattice}); applied to the special case in which
  $w$ is another reflection, this gives (R4).  For (R3) we take the
  following argument from \cite[V.1]{rSerre}. Suppose that $\refl$
  and $\tau$ are reflections in $W$ such that both $b_\refl$ and
  $b_\tau$ are multiples of a single element $b\in L$; we must show
  that $\refl=\tau$. Let $w=\refl\tau$. Then $w(b)=b$, and $w$ acts
  as the identity on $\Q\Tensor(L/\langle b\rangle)$; in particular,
  all of the eigenvalues of $w$ are equal to~$1$. But $w\in W$ has finite
  order, and so $w$ must be the identity map.

  On the other hand, suppose that $R$ is a root system in $L$. For
  each root $ r $ let $\refl_ r \co L\to L$ be given by
  \[
     \refl_r(x) = x + n_r(x)r.
  \]
  It follows from (R2) that $\refl_ r $ is a
  reflection on $L$, with a marking $\pm(r,n_r)$; take $W$ to be the
  group generated the transformations $\refl_r$. By (R4) each
  reflection $\refl_ r $ preserves the finite set $R$ of roots, and
  so $W$ also preserves $R$.  To see that $W$ is finite, it is enough
  to check that any element of $W$ which acts trivially on $R$ is the
  identity, but this follows from (R1) and the fact that the subspace
  $\bigcap_r\ker(n_r)$ is pointwise fixed by each reflection $\refl_r$
  and hence pointwise fixed by $W$. Note that by the above eigenvalue
  argument, $w\refl_rw^{-1}=\refl_{w(r)}$: both transformations send
  $w(r)$ to $-w(r)$ and act as the identity on $\Q\Tensor (L/\langle
  w(r)\rangle)$. Checking that the markings are preserved under
  conjugation thus amounts to showing that if $ \refl_r = \refl_{t}$
  then $r=\pm t$ (note that as in Remark~\ref{HowManyMarkings} a marking
  $\pm(r,n_r)$ is determined by its first component). But since $r$
  and $t$ then determine markings of the same reflection, it follows
  from Remark~\ref{HowManyMarkings} that $r=\pm t$, $r=\pm2t$, or $t=\pm2r$,
  and the last two possibilities are ruled out by (R3).  The one
  remaining issue is that there might be reflections in $W$ which do
  not appear in the set $\{\refl_r\mid r\in R\}$; these reflections
  have not been marked in the above discussion. But by
  \cite[1.14]{rHump} there are no such additional reflections.
\end{proof}

Proposition~\ref{InterpretMarkedTorus} is a consequence of the
following lemma, which is proved by examining the discussion in
Remark~\ref{HowManyMarkings}.

\begin{lem}
  Let $\refl$ be a reflection on the lattice $L$, treated also
  as a reflection on $T=\Torus(1)\Tensor L=(\R\Tensor L)/L$. Then
  sending $b\in L$ to the residue class of $b/2$ in $T$ gives a
  bijection between markings $\pm(b,\beta)$ of $\refl$ and elements
  $h\in T$ satisfying the conditions of Definition~\ref{DTorusMarking}.
\end{lem}

\begin{numbered}{Duality}\label{DualMarkings}
  A reflection $\refl$ on $L$ induces a reflection on the dual
  lattice $L\dual=\Hom(L,\Z)$; sending $\pm(b,\beta)$ to
  $\pm(\beta,b)$ gives a bijection between markings of $\refl$ and
  markings of the dual reflection.  It follows that if $(L,W)$ is a
  marked reflection lattice with markings
  $\pm(b_\refl,\beta_\refl)$, then $(L\dual,W)$ is a marked
  reflection lattice with markings $\pm(\beta_\refl,b_\refl)$.
  Similarly, interchanging roots and coroots give a bijection between
  root systems in $L$ and root systems in $L\dual$.
\end{numbered}

\begin{numbered}{Counting root systems} 
  As in Remark~\ref{HowManyMarkings}, the number of markings of a reflection
  lattice $(L,W)$ is $2^k$, where $k$ is the number of conjugacy
  classes of reflections in $W$ which are trivial mod~$2$. By
  Proposition~\ref{InterpretRootSystem}, this is also the number of distinct root
  systems with reflection lattice $(L,W)$.
\end{numbered}
\end{NumberedSubSection}

\begin{NumberedSubSection}{Another notion of root
    system}\label{TwoRootSystems} We point out the relationship
  between our notion of root system (which is Tits') and that of Serre and Bourbaki.
  Recall that what we term a \emph{geometric root system}
  (Remark~\ref{PointOutDifference}) is a pair $(V,R)$ such that $V$ is a
  finite-dimensional real vector space, $R$ is a finite subset of $V$,
  and the following conditions are satisfied \cite[Section~V]{rSerre}.
\begin{itemize}
\item $R$ generates $V$ as a vector space, and does not contain~$0$.
\item For each $r\in R$ there exists a reflection $s_r$ on
  $V$ which carries $R$ to itself and has $s_r(r)=-r$.
\item For each $r,t\in R$, $s_r(t)-t$ is an
  integral multiple of $r$.
\end{itemize}
The geometric root system is said to be \emph{reduced} if for each
$r\in R$, $r$ and $-r$ are the only roots proportional
to $r$.

The following proposition is immediate. We will say that a root system
(see Definition~\ref{DefineRootSystem}) is \emph{semisimple} if $\bigcap_r\ker(n_r)=\{0\}$. 

\begin{prop}
  If $(L,R,\{n_r\})$ is a semisimple root system, then $(\R\Tensor L,
  R)$ is a a reduced geometric root system.
\end{prop}

Conversely, if $(V,R)$ is a reduced geometric root system, then for
each $r\in R$ define a homomorphism $H_r\co V\to\R$ by the
formula $s_r(x)=x+ H_r(x)r$. Let $L_{\text{max}}$ be
the lattice in $V$ given by $\{x\mid \forall r,
H_r(x)\in\Z\}$, and let $L_{\text{min}}$ be the lattice in $L$
spanned by $R$. Note that if $L\subset V$ is a lattice between
$L_{\text{min}}$ and $L_{\text{max}}$, then $H_r(L)\subset\Z$.

\begin{prop}\label{PIntegralForms}
  Suppose that $(V,R)$ is a reduced geometric root system, and that
  $L$ is a lattice in $V$ with $L_{\text{min}}\subset L\subset
  L_{\text{max}}$. Then $(L,R,\{H_r\})$ is a semisimple root system.
\end{prop}

Recall that a connected compact Lie group $G$ is said to be
\emph{semisimple} if $\pi_1G$ is finite, or equivalently if the center
of $G$ is finite. It is well--known \cite{rBourbaki} that geometric
root systems classify connected semisimple compact Lie groups up to
isogeny, ie up to finite covers. A glance at Serre's description of
the fundamental group of a connected compact Lie group
\cite[VIII-10]{rSerre} shows that extra lattice data in a semisimple
root system (Proposition~\ref{PIntegralForms}) allow such systems to classify
connected semisimple compact Lie groups up to isomorphism.

More generally, root systems classify connected compact Lie groups up
to isomorphism: the lattice $\bigcap_rn_r$ in the root system derived
from a connected compact Lie group $G$ determines the torus which is
the identity component of the center of $G$. In our way of assigning a
root system to $G$ (see Section~\ref{CCompactLie}), which is dual to the
ordinary one (see Remark~\ref{MarkingsAndClasicalRoots}), $\bigcap_rn_r$ is the
fundamental group of the central torus. This classification theorem is
stated in \cite[\S8, Sect.~9]{rBourbakiNine}. It can also
be obtained by combining the correspondence between compact Lie groups
and reductive complex algebraic groups \cite[p247]{rOV} with the
classification of connected reductive complex algebraic groups by
\emph{root data} \cite{rSpringerLAG}. Springer's root data are essential
the same as our root systems.

\end{NumberedSubSection}


\section{The normalizer and reflection extensions}\label{CNormalizerExt}

In this section we take a marked reflection torus $(T,W,\{h_\refl\})$
and construct the associated normalizer extension of $W$ by $T$.  Let
$\refls\subset W$ be the set of reflections in $W$; the group $W$ acts
on $\refls$ by conjugation, and hence on the free abelian group
$\Z[\refls]$ generated by $\refls$. We first describe an extension of
$W$ by $\Z[\refls]$, called the \emph{reflection extension} of $W$,
which depends only on the structure of $W$ as a reflection group, and
then use this in conjunction with the markings $\{h_\refl\}$ to obtain
the normalizer extension.

\begin{NumberedSubSection}{The reflection extension}\label{DReflectionExtension}
  Write $\refls=\coprod_i \refls_i$ as a union of conjugacy classes of
  reflections. For each index $i$ choose a reflection $\crefl_i$ in
  $\refls_i$, and let $C_i$ be the centralizer of $\crefl_i$ in $W$.
  Let $L=\pi_1T$, and let $a_i\in L$ be a nonzero element with
  $\crefl_i(a_i)=-a_i$. It is clear that each element of $C_i$ takes
  $a_i$ to $\pm a_i$, and that $C_i$ is isomorphic to $\langle
  \crefl_i\rangle \times C_i\pperp$, where $C_i\pperp$ is the subgroup
  of $C_i$ consisting of elements which fix $a_i$.

As usual,  extensions of $W$ by $\Z[\refls]$ are classified by elements of
$\HH^2(W;\Z[\refls])$. There are isomorphisms
  \[
   \HH^2(W;\Z[\refls])\cong \oplus_i\HH^2(W;\Z[\refls_i])\cong\oplus_i\HH^2(C_i;\Z)
  \]
  where the last isomorphism comes from Shapiro's lemma. Thus in order
  to specify an extension of $W$ by $\Z[\refls]$ it is enough to specify
  for each $i$ an extension of $C_i$ by $\Z$. We do this as follows.
  For each $i$ write
  \[
     C_i\cong \langle \crefl_i\rangle \times C_i\pperp \cong \Z/2\times C_i\pperp
  \]
  and consider the obvious extension
  \begin{equation}\label{KeyLocalExtension}
  1 \to \Z \to
  \Z\times C_i\pperp \to \Z/2\times C_i\pperp \to 1.
  \end{equation}
  This extension corresponds to the element of $\HH^2(C_i;\Z)$ which
  is a pullback under the projection $C_i\to\Z/2$ of the unique
  nonzero element in $\HH^2(\Z/2;\Z)$. By Shapiro's lemma there is a
  corresponding extension $\RE_i(W)$ of $W$ by $\Z[\refls_i]$, and it
  is easy to see that up to canonical isomorphism this extension of
  $W$ does not depend on the choice of representative $\crefl_i$ for
  the conjugacy class $\refls_i$.

  \begin{defn}\label{DefineReflectionExtension}
    The \emph{reflection extension} $\RE( W)$ of $W$ is the extension
    of $W$ by $\Z[\refls]$ determined by the above collection of extensions
    of the groups $C_i$; in other words, $\RE(W)$ is the fibrewise
    product (over $W$) of the
    extensions
    \[
         1 \to \Z[\refls_i] \to \RE_i(W) \to W \to 1.
    \]
  \end{defn}

\end{NumberedSubSection}

\begin{NumberedSubSection}{The normalizer extension}
  Given a marked reflection torus
$$\mcT=(T,W,\{h_\refl\}),$$
there is a $W$--map
\[
    f^{\mathcal T}\co  \Z[\refls]\to T
\]
which sends a reflection $\refl$ to $h_\refl$.

\begin{defn}\label{DefineNormalizerExt}
  Suppose that $\mcT=(T,W,\{h_\refl\})$ is a marked reflection torus.
  The \emph{normalizer extension} $\Norm(T,W,\{h_\refl\})$ of $W$ by
  $T$ is defined to be the image of the reflection extension $\RE(W)$
  under the map $f^{\mcT}\co \Z[\refls]\to T$.
\end{defn}

\begin{rem}
  The reflection extension of $W$ is determined by an extension class
  $k\in \HH^2(W;\Z[\refls])$, and the normalizer extension by
  $f^{\mcT}_*(k)\in \HH^2(W;T)$.
\end{rem}

\end{NumberedSubSection}

\section{The Tits extension}\label{CExtended}

Suppose that $(L,W)$ is a reflection lattice.  In this section we
derive generators and relations (Section~\ref{GeneratorsAndRelations}) for the
reflection extension $\RE(W)$ of $W$. To accomplish this, we identify
$\RE(W)$ with the \emph{extended Weyl group} constructed by Tits
\cite{rTits}; we call this extended Weyl group the \emph{Tits
  extension} of $W$ and denote it $\Tits(W)$. The presentation for
$\RE(W)$ is used both in the proof of
Theorem~\ref{IdentifyTorusNormalizer} and later on in the proof of
Proposition~\ref{IdentifyNT} (Coxeter case).

\begin{NumberedSubSection}{The Tits extension}
  Let $\srefls=\{\srefl_1,\ldots,\srefl_\ell\}$ be a chosen set of simple
  reflections in $W$ (see Remark~\ref{CoxeterGeneration}), and let $\refls^*$ be
  a copy of the set $\refls$ of all reflections in $W$. As in
  Section~\ref{CNormalizerExt}, $W$ acts by conjugation on $\refls$ and thus
  on $\Z[\refls]$ or $\Z[\refls^*]$.

\begin{defn}
  The \emph{Tits extension} $\Tits(W)$ of $W$ is the subgroup of 
  the semidirect product $\Z[\refls^*]\rtimes
  W$ generated by the elements $(\srefl^*,\srefl)$, $\srefl\in \srefls$.
\end{defn}

Since the set $\srefls$ generates $W$ \cite[1.5]{rHump}, the natural
surjection $\Z[\refls^*]\rtimes W\to W$ restricts to a surjection
$\Tits(W)\to W$.  The kernel of this map is contained in
$\Z[\refls^*]$, and Tits \cite[Cor.~2.7]{rTits} shows that it is exactly $2\Z[\refls^*]$,
which of course is isomorphic to $\Z[\refls]$ as a $W$--module. Under
this identification, we can treat $\Tits(W)$ as an extension
\begin{equation}\label{TitsExtension}
  1 \to \Z[\refls]\to \Tits(W)\to W \to 1
\end{equation}
of $W$ by $\Z[\refls]$. In this section we will prove the following
proposition.

\begin{prop}\label{TwoExtensionsMatch}\label{IdentityTheorem}
  Suppose that $(L,W)$ is a reflection lattice and that $\refls$ is the set
  of reflections in $W$. Then, as an extension of $W$ by $\Z[\refls]$, the
  reflection extension $\RE(W)$ is isomorphic to the Tits extension
  $\Tits(W)$.
\end{prop}

\begin{rem}
  It follows from 
  Proposition~\ref{TwoExtensionsMatch} that  the isomorphism class of the
  extension $\Tits(W)$ is independent of the choice of the set
  $\srefls$ of simple reflections. This also follows from the fact
  that any two sets of simple reflections are conjugate in~$W$
  \cite[1.4]{rHump}. 
\end{rem}

\begin{numbered}{Generators and relations}\label{GeneratorsAndRelations}
  Let $(m_{i,j})$ be the Coxeter matrix (see Remark~\ref{CoxeterGeneration}) of
  $(L,W)$, so that $W$ is generated by simple reflections
  $\{\srefl_1,\ldots,\srefl_\ell\}$ subject to the relations
  \begin{equation}\label{CoxeterRelations}
  (\srefl_i\srefl_j)^{m_{i,j}}=1.
  \end{equation}
  Note that $m_{i,i}=1$. Following Tits \cite[0.1]{rTits}, we put
  these relations in a slightly different form.  Write
  \[\prodd(n;y,x) = \cdots yxyx\]
  where the product on the right contains exactly $n$ 
  factors, alternating between $x$ and $y$, and beginning on the right with an $x$. Then $W$
  is generated by the elements $\srefl_i$, $1\le i\le \ell$, subject only
  to the relations
  \[
  \begin{aligned}
  \srefl_i^2&=1\\
  \prodd(m_{i,j};\srefl_i,\srefl_j)&=\prodd(m_{i,j};\srefl_j,\srefl_i) \text{ for }i\ne
  j.
  \end{aligned}
  \]
  In fact the first relation gives $\srefl_i=\srefl_i^{-1}$, while the
  second one is obtained from \eqref{CoxeterRelations} by the reversible
  process of moving the leading terms of the expanded expression
  $(\srefl_i\srefl_j)^{m_{i,j}}$ to the right hand side of the
  equation and replacing each occurrence of $x^{-1}$ ($x= \srefl_i$
  or~$\srefl_j$) by $x$ itself.

  In \cite[2.6]{rTits}, Tits gives a presentation for $\Tits(W)$
  parallel to the above presentation for $W$. For notational reasons,
  think for the moment of $\Tits(W)$ as an extension of $W$ by
  $\Z[\refls^*]$.  Then $\Tits(W)$ is generated by $\Z[\refls^*]$ and
  symbols $q_i$, one for each simple reflection, subject to the
  following relations:
\[
\begin{aligned}
    q_i^2 &= \srefl^*_i\\
    q_i\refl^*q_i^{-1}&=(\srefl_i\refl\srefl_i^{-1})^*, \quad \refl\in\refls\\
    \prodd(m_{i,j};q_i,q_j)&=\prodd(m_{i,j};q_j,q_i)\text{ for }i\ne
    j.
\end{aligned}
\]
With respect to this presentation of $\Tits(W)$, the quotient map
$\Tits(W)\to W$ sends $q_i$ to $\srefl_i$. By
Proposition~\ref{TwoExtensionsMatch}, these same formulas give a presentation of
$\RE(W)$.
\end{numbered}
  
\end{NumberedSubSection}

The rest of this section is devoted to the proof of Proposition~\ref{TwoExtensionsMatch}.
We first need an explicit formulation of Shapiro's lemma.
As above, if $G$ is a finite group and $H\subset G$ is a subgroup,
extensions of $G$ by the permutation module $\Z[G/H]$ correspond
bijectively to extensions of $H$ by the trivial $H$--module $\Z$. We
will call these two extensions \emph{Shapiro companions} of one
another. The following lemma is elementary; see \cite[III.8, Ex.~2]{rBrown}.

\begin{lem}\label{ExplicitShapiro}
  Suppose that $G$ is a finite group, $H\subset G$ a subgroup, and 
  \begin{equation}\label{BigExtension}
    1 \to \Z[G/H] \to E \to G \to 1
  \end{equation}
  an extension of $G$ by the permutation module $\Z[G/H]$. Then the
  Shapiro companion of \eqref{BigExtension} is obtained by pulling the
  extension back over $H$ to obtain
  \[
    1 \to \Z[G/H] \to E' \to H \to 1,
  \]
 and then taking the quotient of $E'$ by the subgroup of
  $\Z[G/H]$ generated by the elements $\{x\in G/H, x\ne eH \}$.
\end{lem}

\begin{rem}
  Note that $H$ fixes the coset $eH$ in $G/H$, and so $H$ carries the set
  $\{x\in G/H, x\ne eH\}$ to itself. This implies that the subgroup of $\Z[G/H]$
  generated by this set is a normal subgroup of the group $E'$ above.
  The quotient
  \[ \Z[G/H]\,\,/\,\,\Z[x\in G/H, x\ne
  eH]
  \]
  is canonically isomorphic to the trivial $H$--module $\Z$, with
  generator given by the residue class of $eH$.
\end{rem}

We next recall some apparatus from \cite{rTits}. Let \Iset/ be a set
of indices corresponding to the chosen simple reflections in $W$,
and \Igroup/ the free group on \Iset/.  An element  $\Iel\in\Igroup$
is said to be \emph{positive} if it is a product of elements of \Iset;
the \emph{length} $l(\Iel)$ is then the number of elements in this
product expression. To conform with \cite{rTits}, for each $i\in\Iset$
we will denote the corresponding simple reflection in $W$ by $r_i$
(instead of $\srefl_i$); there is a unique (surjective) homomorphism
$r\co \Igroup\to W$ with $r(i)=r_i$, $i\in \Iset$. A positive word
$\Iel\in\Igroup$ is said to be \emph{minimal} if there does not exist
a positive word $\Jel\in\Igroup$ with $l(\Jel)<l(\Iel)$ and
$r(\Jel)=r(\Iel)$. If $w\in W$, the \emph{length} of $w$ is defined to
be the length of a minimal positive word in \Igroup/ which maps to $w$
under $r$.

Suppose that $\Iel=i_1i_2\cdots i_m$ is a positive word in $\Igroup$
of length $m$. For each integer $k$ with $1\le k\le m$ set
\[
   \refl_k= r(i_1i_2\cdots i_{k-1}) \cdot r_{i_k} \cdot r(i_1i_2\cdots i_{k-1})^{-1}.
\]
Each one of these elements is a reflection in $W$, and the sequence
$\refl_1, \refl_2,\cdots ,\refl_m$ is said to be the \emph{sequence of
  reflections associated to the word $\Iel$}. Recall the following
proposition.

\begin{prop}
  {\rm\cite[1.6]{rTits}}\qua Let $\Iel\in\Igroup$ be a minimal word, $\refl\in W$ a
  reflection, and $n$ the number of times which $\refl$ appears in the
  sequence of reflections associated to $\Iel$. Then $n=0$ or $1$.
\end{prop}

We need to evaluate the integer $n$ above in a special case.

\begin{prop}\label{ZeroProposition}
  Let $\Iel\in\Igroup$ be a minimal word, $\refl\in W$ a
  reflection, and $n$ the number of times which $\refl$ appears in the
  sequence of reflections associated to $\Iel$. Suppose that $\refl$ is a
  simple reflection, that $r(\Iel)$ is a reflection which
  commutes with $\refl$, and that $r(\Iel)\ne \refl$. Then $n=0$.
\end{prop}

The proof of this was explained to us by M. Dyer. It depends on the
following lemma.

\begin{lem}\label{DyerLemma}
  {\rm\cite[1.4]{rDyer}}\qua Suppose that
  $t=r(i_1i_2\cdots i_{2n+1})$ is a reflection in $W$ with
  $l(t)=2n+1$. Then $t=r(i_1\cdots i_ni_{n+1} i_n\cdots i_1)$.
\end{lem}

For the convenience of the reader, we give the argument from
\cite{rDyer}.

\begin{proof}[Proof of Lemma~\ref{DyerLemma}]
  Let $x=r(i_n\cdots i_1)$ and $y=r(i_{n+2}\cdots i_{2n+1})$. Then
  $l(x)=l(y)=n$ and $l(r_{i_{n+1}}x)=l(r_{i_{n+1}}y)=n+1$.  Note that
  $t=t^{-1}$, so that $r_{i_{n+1}}yt=r_{i_{n+1}}yt^{-1}=x$ and
  $l(r_{i_{n+1}}yt)<l(r_{i_{n+1}}y)$.  Applying the Strong Exchange
  Condition \cite[5.8]{rHump} to this last observation shows that
  there is a number $k$ between $n+1$ and $2n+1$ such that
\[
    x = r_{i_{n+1}}yt= r(i_{n+1}\cdots \widehat {i_k}\cdots i_{2n+1}).       
\]
Since $l(x)=n$, this is a minimal expression for $x$. Since
$l(r_{i_{n+1}}x)>l(x)$ it must be the case that $k=n+1$ and hence that
$x=y$.
\end{proof}

\begin{proof}[Proof of Proposition~\ref{ZeroProposition}]
Write $\Iel=i_1i_2\cdots i_m$, and suppose that $\refl$ appears among the
sequence $\refl_1, \refl_2,\cdots, \refl_m$ of reflections associated to \Iel,
that is, suppose that there exists $k$ with $1\le k\le m$ such
that
\[
  \refl= \refl_k = r(i_1\cdots i_{k-1}) r_{i_k} r(i_1\cdots i_{k-1})^{-1}.
\]
Let $t=r(i_1\cdots i_m)$; it is clear that
\[
\begin{aligned}
   t &= r(i_1\cdots i_{k-1}) \cdot r_{i_k} \cdot r(i_1\cdots i_{k-1})^{-1} 
        r(i_1\cdots i_{k-1})r(i_{k+1}\cdots i_m)\\
     & = \refl \cdot r(i_1\cdots i_{k-1})r(i_{k+1}\cdots i_m).
\end{aligned}
\]
Choose $i_0\in \Iset$ so that $\refl=r_{i_0}$ (recall that $\refl$ is  a
simple reflection) and let $\Jel = i_0i_1\cdots
i_{k-1}i_{k+1}\cdots i_m$; then $r(\Jel)=t$ and, since
$l(\Jel)=l(\Iel)$, $\Jel$~is a minimal word with $r(\Jel)=t$. Write
$\Jel= j_1\cdots j_m$, with $j_1=i_0$, so that $r_{j_1}=\refl$. For
determinant reasons $m=2q+1$ is odd, and $m>1$ because $t=r(\Jel)\ne
\refl$. By Lemma~\ref{DyerLemma},
\[
\begin{aligned}
    t&= r_{j_1}\cdot r(j_2\cdots j_qj_{q+1}j_q\cdots j_2) \cdot r_{j_1}\\
     &= \refl \cdot r(j_2\cdots j_qj_{q+1}j_q\cdots j_2)\cdot  \refl.
\end{aligned}
\]
Now consider $\refl\cdot t\cdot \refl$. Since $t$ commutes with $\refl$ and
$\refl^2=1$, this element is equal to $t$. In combination with the above
equation, we get
\[
\begin{aligned}
    t &= \refl\cdot t\cdot \refl\\
      &= \refl \cdot \refl\cdot r(j_2\cdots j_qj_{q+1}j_q\cdots j_2)\cdot
      \refl\cdot \refl\\
      & = r(j_2\cdots j_qj_{q+1}j_q\cdots j_2),
\end{aligned}
\]
which, in contradiction to the assumption that $\Iel$ is a minimal
word, gives $l(t)\le m-2$. This contradiction shows that $\refl$ does not
appear in the list of reflections associated to $\Iel$.
\end{proof}

Finally, one more calculation from Tits. Let
$a\co \Igroup\to\Z[S^*]\rtimes W$ be the homomorphism determined by
$a(i)= (r_i^*,r_i)$, $i\in\Iset$. By definition, the image of $a$ is
$\Tits(W)$. 

\begin{lem}{\rm\cite[Lemme~2.4]{rTits}}\qua\label{TitsCalculation}
  Let $\Iel\in\Igroup$ be a positive word of length $m$, and
  $\refl_1,\refl_2,\cdots,\refl_m$ the sequence of reflections associated to
  \Iel. Then 
  \[
     a(\Iel)=\left(\sum_{k=1}^m \refl_k^* , r(\Iel)\right).
  \]
\end{lem}

\begin{proof}[Proof of Proposition~\ref{IdentityTheorem}]
  We identify the kernel of
  the map $\Tits(W)\to W$ with $\Z[\refls]$. Then $\Tits(W)$ sits inside
  the semidirect product $\Z[\refls^*]\rtimes W$, and there is map of short
  exact sequences
  \[
  \begin{CD}
    1 @>>> \Z[\refls] @>>>\Tits(W) @>>>W@>>>1\\
@. @VVV @VVV @VVV @.\\
    1 @>>> \Z[\refls^*] @>>>\Z[\refls^*]\rtimes W@>>>W@>>>1
  \end{CD}
  \]
  in which the left hand vertical map gives an isomorphism between
  $\Z[\refls]$ and $2\Z[\refls^*]$. As in Section~\ref{DReflectionExtension},
  giving the extension $\Tits(W)$ amounts to giving a collection of
  extensions $\{\Tits_i(W)\}$, where $\Tits_i(W)$ is an extension of
  $W$ by $\Z[\refls_i]$; the extension $\Tits_i(W)$ is obtained from
  $\Tits(W)$ by taking a quotient of $\Tits(W)$ by the subgroup of
  $\Z[\refls]$ generated by $\refls\setminus \refls_i$. Clearly
  $\Tits_i(W)$ lies inside the semidirect product $\Z[\refls_i^*]\rtimes
  W$, and in fact there is a map of short exact sequences
 \begin{equation}\label{CommDiag}
  \begin{CD}
    1 @>>> \Z[\refls_i] @>>>\Tits_i(W) @>>>W@>>>1\\
@. @VVV @VVV @VVV @.\\
    1 @>>> \Z[\refls_i^*] @>>>\Z[\refls_i^*]\rtimes W@>>>W@>>>1
  \end{CD}
  \end{equation}
  in which as before the left hand vertical map identifies $\Z[\refls_i]$
  with $2\Z[\refls_i^*]$. 
  Pulling
  the diagram \eqref{CommDiag} back over the inclusion $C_i\to W$
  gives a commutative diagram:
\begin{equation}\label{Ei}
  \begin{CD}
    1 @>>> \Z[\refls_i] @>>>E_i @>>>C_i@>>>1\\
@. @VVV @VVV @VVV @.\\
    1 @>>> \Z[\refls_i^*] @>>>\Z[\refls_i^*]\rtimes C_i@>>>C_i@>>>1
  \end{CD}
\end{equation}
 Now taking quotients as in Lemma~\ref{ExplicitShapiro}, in the upper group
 by $\Z[\refls_i\setminus \crefl_i]$ and in the lower by $\Z[\refls_i^*\setminus
 \crefl_i^*]$, yields:
\begin{equation}\label{Eiprime}
  \begin{CD}
    1 @>>> \Z \crefl_i @>>>E_i' @>>>C_i@>>>1\\
@. @VVV @VVV @VVV @.\\
    1 @>>> \Z \crefl_i^* @>>>\Z \crefl_i^*\times C_i@>>>C_i@>>>1
  \end{CD}
\end{equation}
It is clear that the left hand vertical map sends $\crefl_i$ to $2\crefl_i^*$.
We need to argue that the upper extension is the one described in
\eqref{KeyLocalExtension}. 

Note that every reflection in $W$ is conjugate to a simple
reflection \cite[1.14]{rHump}, so we can assume that $\crefl_i$
belongs to the set $\srefls$ of simple reflections.  Recall that
$C_i\cong \langle \crefl_i\rangle \times C_i^\perp$, where $C_i^\perp$
can be identified as the subgroup of $W$ which pointwise fixes
$L^-(\crefl_i)$. Following
Steinberg \cite[1.5]{rSteinbergDifferential} \cite[1.12]{rHump},
$C_i^\perp$ is generated by the reflections it contains, and so in
particular it is generated by reflections $t$ which commute with $\crefl_i$
and are distinct from $\crefl_i$.  Pick such a $t$, and write
$t=r(\Iel)$, where $\Iel$ is a minimal positive word. By
Proposition~\ref{ZeroProposition}, $\crefl_i$ does not appear among the list of
reflections associated to $\Iel$. Clearly the image $\alpha$
of $a(\Iel)$ in $\tau_i(W)$ belongs to $E_i$ (see
\eqref{Ei}); let $\bar a(\Iel)$ be the image of $\alpha$ in $E_i'$ (see
\eqref{Eiprime}).  It now follows from Lemma~\ref{TitsCalculation}
that the image of $\bar a(\Iel)$ in $\Z \crefl_i^*\times C_i$ is
$(0,t)$. The subgroup of $\Z \crefl_i^*\times C_i$ generated by these
elements as $t$ varies through the reflections in $C_i^\perp$ is
$\{0\}\times C_i^\perp$. Consider the subgroup $E_i''$ of $\Z
\crefl_i^*\times C_i$ generated by $\{0\}\times C_i^\perp$ and by
$(\crefl_i^*,\crefl_i)$. Clearly, $E_i''$ is contained inside the
image of $E_i'$ in $ \Z \crefl_i^*\times C_i$; on the other hand, it
is also easy to check that there is a short exact sequence
\[
 1\ \to 2\Z \crefl_i^* \to E_i'' \to C_i \to 1.
\]
By the five lemma, $E_i''\cong E_i'$, and it is now simple to verify that
$E_i''$ is the extension of $C_i$ described in
\eqref{KeyLocalExtension}. 
\end{proof}

\section{Compact Lie Groups}\label{CCompactLie}
In this section, we follow the lead of Tits \cite{rTits} and describe
the normalizer of the torus in a connected compact Lie group in terms
of a marked reflection torus (Definition~\ref{DTorusMarking}) associated to the
group.

Suppose that $G$ is a connected compact Lie group of rank~$r$ with
maximal torus $T$, maximal torus normalizer $N=NT$ and Weyl group
$W=N/T$. Then $T\iso \Torus(r)$ and the conjugation action of $W$ on
$T$ identifies $W$ with a subgroup of $\Aut(T)$ generated by
reflections \cite[VII.2.13]{rHelgason} \cite[5.16]{rDWelementary}. In
particular, $(T,W)$ is a reflection torus.

For each reflection $\refl\in W$ we will construct a
subgroup $N_\refl$ of $N$, called the \emph{root subgroup} of $N$
associated to $\refl$. We will also single out a specific element
$h_\refl$ in $T$, the \emph{marking} associated to $\refl$.

Given $\refl$, let $T^+_0(\refl)$ denote the identity component of the
fixed point set of the action of $\refl$ on $T$. It is easy to see
that $T^+_0(\refl)$ is a torus of rank $r-1$; if $L$ is the lattice
corresponding to $T$, then under the isomorphism $\Torus(1)\Tensor
L\iso T$, $T^+_0(\refl)$ is given by $\Torus(1)\Tensor L^+(\refl)$.
Let $\RG\refl$ denote the centralizer of $T^+_0(\refl)$ in $G$. Then
$\RG\refl$ is a connected \cite[VII.2.8]{rHelgason}
\cite[7.3]{rDWelementary} compact Lie subgroup of $G$ with maximal
torus $T$ and Weyl group given by the subgroup $\{1,\refl\}$ of $W$.
The normalizer of $T$ in $\RG\refl$ can be identified with the
centralizer of $T^+_0(\refl)$ in $N$; we denote this group by
$N(\refl)$.  The only simply-connected compact Lie group of rank~$1$
is $\SU(2)$; the center of $\SU(2)$ is $\Z/2$ and the adjoint form is
$\SO(3)$.  Given this, it follows from the product splitting theorem
for compact Lie groups \cite[Section~10]{rDWelementary} that $\RG\refl$ is
isomorphic to one of the following three groups:
\begin{enumerate}
\item A product $\SU(2)\times \Torus(r-1)$
\item A product $\SO(3)\times \Torus(r-1)$
\item A product $\operatorname{U}(2)\times \Torus(r-2)$
\end{enumerate}

Recall that $T^-_0(\refl)$ denotes the identity component of
$T^-(\refl)$; this group is a torus of rank~$1$. 

\begin{defn}\label{DefineLGRootSubgroup}
 The \emph{root
subgroup} $N_\refl$ is given by
\[
   N_\refl=\{x\in N(\refl)\mid x\text{ is conjugate in }\RG\refl\text{
     to some }y\in T^-_0(\refl)\}.
\]
\end{defn}

\begin{lem}\label{RootSubgroupIsSubgroup}
  The set $N_\refl$ is a closed subgroup of $N(\refl)$ with two
  components. The identity component of $N_\refl$ is the circle
  $T^-_0(\refl) = N_\refl\cap T$, and the natural map
  $N_\refl/T^-_0(\refl)\to N(\refl)/T\iso \Z/2$ is an isomorphism. 
  If $x, y$ belong to the nonidentity component of $N_\refl$, then
  $x^2=y^2\in T^-_0(\refl)$.
\end{lem}

\begin{proof}
  The simplest way to prove this is to examine the above three
  cases. In the first two cases $N_\refl$ is the normalizer of the
  maximal torus $T^-_0(\refl)$ in $\SU(2)$ (respectively, $\SO(3)$). In the
  third case, $N_\refl$ is the normalizer of the maximal torus
  $T^-_0(\refl)$ of the subgroup $\SU(2)$ of $\operatorname{U}(2)$.
\end{proof}

\begin{defn}\label{DefineLGMarking}
  The \emph{marking} $h_\refl\in T^-_0(\refl)$ is the image of the
  nonidentity component of $N_\refl$ under the squaring map $x\mapsto
  x^2$.
\end{defn}

\begin{lem}\label{LGMarkings}
  The element $h_\refl\in T^-_0(\refl)\subset T$ is
  a marking for $\refl$ in the sense of Definition~\ref{DTorusMarking}.
\end{lem}

\begin{proof}
  Again, this is most easily proved by inspection.  The element
  $h_\refl$ is strongly $\refl$--negative by construction, and it is
  easy to see that $2h_\refl=0$ (where as usual we write the group
  operation in the torus additively).  Finally, if $\refl$ is
  nontrivial mod~$2$, then $\RG\refl\iso
  \operatorname{U}(2)\times\Torus(r-2)$ and $h_\refl$ is the
  nontrivial central element of order~$2$ in $\operatorname{U}(2)$.
\end{proof}

The following statement is a consequence of the naturality of the
above constructions.

\begin{prop}\label{LGConjugationInvariant}
  Suppose that $G$ is a connected compact Lie group with maximal torus
  $T$, maximal torus normalizer $N$, Weyl group $W$, root subgroups
  $\{N_\refl\}$, and markings $\{h_\refl\}$ as above. Suppose that
  $x\in N$ has image $w\in W$. Then for any reflection $\refl\in W$,
  \[
    \begin{aligned}
      xN_{\refl\vphantom{w^{-1}}} x^{-1}&=N_{w\refl w^{-1}}\\
      x h_{\refl\vphantom{w^{-1}}} x^{-1} &= h_{w\refl w^{-1}}
     \end{aligned}
  \]
  In particular, $(T, W,\{h_\refl\})$ is a marked reflection torus. 
\end{prop}

\begin{rem}\label{MarkingsAndClasicalRoots}
  Let $L$ denote $\pi_1T$ and let $R$ be the root system in $L$
   associated
  to the above marked reflection torus (see
  Propositions~\ref{InterpretRootSystem} and
  \ref{InterpretMarkedTorus}). It is not hard to identify the
  dual root system $R\dual$ in $L\dual$ (see Section~\ref{DualMarkings}) with
  the root system usually associated to $G$ \cite[VIII-8, VI]{rSerre}.
  The correspondence $G\mapsto R\dual$ produces a bijection between
  isomorphism classes of connected compact Lie groups and isomorphism
  classes of root systems (see Section~\ref{TwoRootSystems}).

  The marking $h_\refl$ is defined above (Definition~\ref{DefineLGMarking}) in
  terms of $N_\refl$. However, it is not hard to see that $h_\refl$ is
  the only marking for $\refl$ in the intersection of $T^-_0(\refl)$ with the
  square of the nonidentity component of $N(\refl)$. Since $N(\refl)$
  is the centralizer in $N$ of $T^+_0(\refl)$, the element $h_\refl$,
  and hence the marked reflection torus $(T,W,\{h_\refl\})$ and the
  root system $R^\#$, are determined by the group structure of $N$. In
  this way we recover the result from \cite{rCurtisEtAl} and \cite{rOsse} that
  $G$ is determined up to isomorphism by $N$.
\end{rem}

The main result of this section is
the following one.

\begin{thm}\label{IdentifyTorusNormalizer}
  Let $G$ be a connected compact Lie group, and $(T,W,\{h_{\refl}\})$
  the marked reflection torus derived from $G$
  (Proposition~\ref{LGConjugationInvariant}).  Then the normalizer of $T$ in
  $G$ is isomorphic to the normalizer extension
  $\Norm(T,W,\{h_{\refl}\})$ (see Section~\ref{CNormalizerExt}).
\end{thm}

This is a consequence of the following theorem, due essentially to
Tits. See Section~\ref{GeneratorsAndRelations} for the notation ``$\prodd(n;x,y)$''.

\begin{thm}\label{MainTitsTheorem}
  Suppose that $G$ is a connected compact Lie group with maximal torus
  $T$, maximal torus normalizer $N$, Weyl group $W$, root subgroups
  $\{N_\refl\}$, and markings $\{h_\refl\}$.  Choose a set
  $\{\srefl_i\}$  of simple reflections in $W$, and let
  $\{m_{i,j}\}$ be the corresponding Coxeter matrix. For
  each $i$, choose an element $x_i$ in the
  nonidentity component of $N_{\srefl_i}$. Then the elements $x_i$
  satisfy the following relations:
  \begin{equation}\label{XRelations}
    \begin{aligned}
           x_i^2&= h_{\srefl_i}\\
           x_itx_i^{-1} &= \srefl_i(t), \quad t\in T\\
           \prodd(m_{i,j};x_i,x_j)&=\prodd(m_{i,j};x_j,x_i)\text{ for
             }i\ne j.
    \end{aligned}
   \end{equation}
\end{thm}

\begin{proof}[Proof of Theorem~\ref{IdentifyTorusNormalizer} (given
Theorem~\ref{MainTitsTheorem})]
  Let $\Norm(W)$ denote the normalizer extension in question.
  According to Definition~\ref{DefineNormalizerExt}, Proposition~\ref{IdentityTheorem}, and
  Section~\ref{GeneratorsAndRelations}, $\Norm(W)$ is generated by the torus
  $T$ together with symbols $q_i$, one for each simple reflection
  $\srefl_i$, subject to the following relations:
\begin{equation}\label{FinalNTRelations}
\begin{aligned}
      q_i^2 &= h_{\srefl_i}\\
      q_i tq_i^{-1}&= \srefl_i(t)\quad t\in T\\
      \prodd(m_{i,j};q_i,q_j)&=\prodd(m_{i,j};q_j,q_i)\text{ for }i\ne
        j.
\end{aligned}
\end{equation}
By Theorem~\ref{MainTitsTheorem}, there is a homomorphism $\phi\co \Norm(W)\to
N$ which is the identity on $T$ and sends $q_i$ to $x_i$. This is an
isomorphism because it is part of a map
\[
\begin{CD}
  1 @>>> T @>>> \Norm(W) @>>> W @>>> 1\\
  @.     @V\|VV      @V \phi VV   @V\|VV   @.\\
  1 @>>> T @>>> N @>>> W @>>> 1
\end{CD}
\]
of exact sequences.
\end{proof}

The proof of Theorem~\ref{MainTitsTheorem} depends on two auxiliary results.
Suppose that $\tilde G\to G$ is a finite covering of connected compact
Lie groups, and that $\tilde T$ is a maximal torus in $\tilde G$. Then
the image $T$ of $\tilde T$ in $G$ is a maximal torus for $G$, the
normalizer $\tilde N$ of $\tilde T$ in $\tilde G$ maps onto the
normalizer $N$ of $T$ in $G$, and the induced map $\tilde N/\tilde
T\to N/T$ is an isomorphism.  In particular, $G$ and $\tilde G$ have
isomorphic Weyl groups, which we will denote by the same letter $W$.
The lattice corresponding to $\tilde T$ is a subgroup of finite index
in the lattice corresponding to $T$, and it follows that an element of
$W$ acts as a reflection on $\tilde T$ if and only if it acts as a
reflection on $T$.  The following statement is clear.

\begin{prop}\label{CoveringStuff}
  Suppose that $p\co \tilde G\to G$ is a finite covering of connected
  compact Lie groups with Weyl group $W$. Let $\tilde T$ be a maximal
  torus in $\tilde G$ and $T\subset G$ the image maximal torus in $G$.
  For each reflection $\refl\in W$, let $\tilde N_\refl$, $\tilde
  h_\refl$ (respectively, $N_\refl$, $h_\refl$) denote the root subgroup
  and marking in $\tilde G$ (respectively, $G$) corresponding to $\refl$.
  Then $p(\tilde N_\refl)=N_\refl$ and $p(\tilde
  h_\refl)=h_\refl$.
\end{prop}

\begin{prop}\label{IntersectionLemma}
  Suppose that $G$ is a connected, simply connected compact Lie group
  with maximal torus $T$, torus normalizer $N$, Weyl group $W$, root
  subgroups $\{N_\refl\}$, and markings $\{h_\refl\}$.  Let
  $\{\srefl_1,\cdots,\srefl_{\ell}\}$ be a set of simple
  reflections for $W$. Then if $i$ and $j$ are distinct integers
  between $1$ and $\ell$, $N_{\srefl_i}\cap N_{\srefl_j}=\{1\}$.
\end{prop}

\begin{proof}
  The image of $N_\refl$ in $W$ is $\{1,\refl\}$, and so if $i\ne j$
  the intersection $N_{\srefl_i}\cap N_{\srefl _j}$ is contained in
  $T$.  Since $N_{\refl}\cap T$ is the identity component
  $T^-_0(\refl)$ of $T^-(\refl)$, we must show that for $i\ne j$,
  $T^-_0(\srefl_i)\cap T^-_0(\srefl_j)=\{0\}$.

  Consider the root system $R$ in $L=\pi_1T$ associated to the marked
  reflection torus $(T,W,\{h_\refl\})$
  (Remark~\ref{MarkingsAndClasicalRoots}). It follows from the fact that $G$
  is simply connected that the roots in $R$ (equivalently, the coroots
  in $R\dual$) span $L$; see \cite[VIII.1]{rSerre}, where $\Gamma(G)$
  denotes what we call $L$ and $\{H_\alpha\}$ denotes the set of
  coroots in $R^\#$. Let $a_i$ be the root corresponding to
  $\srefl_i$. Since all of the roots in $R$ are integral linear
  combinations of the $a_i$ \cite[V.8]{rSerre}, it is clear that the
  $a_i$ span $L$. In view of the fact that $G$ is simply connected,
  the number of simple roots is equal to the rank of $G$, and so
  $\{a_1,\ldots,a_\ell\}$ form a basis for $L$.  It follows
  immediately that $L^-(\srefl_i)=\langle a_i\rangle$, and so
  $L\iso\oplus_i L^-(\srefl_i)$. Tensoring with $\Torus(1)$ gives that
  $T\iso\oplus_i T_0^-(\refl_i)$, which gives the desired result.
\end{proof}

 \begin{proof}[Proof of Theorem~\ref{MainTitsTheorem}]
   If the theorem is true for $G$, it is  true for
   the product $G\times\Torus(n)$. By Proposition~\ref{CoveringStuff}, if the
   theorem is true for $G$ and $A$ is a finite subgroup in the center
   of $G$, it is also true for $G/A$. Now any connected
   compact Lie group $G$ can be written as $(K\times\Torus(n))/A$,
   where $K$ is simply connected and $A$ is a finite subgroup of the
   center of $K\times \Torus(n)$. It is enough, then, to work in the
   case in which $G$ is simply connected.

   In $\SU(2)$, the square of any element in the nonidentity component
   of the torus normalizer is equal to the unique nontrivial central
   element of the group, so it follows from the choice of $h_{\srefl_i}$ that
   $x_i^2=h_{\srefl_i}$. Since $x_i\in N$ projects to $\srefl_i\in W$, it is also
   clear that for $t\in T$, $x_itx_i^{-1}=\srefl_i(t)$. What remains to
   check is whether
    \begin{equation}\label{FundamentalRelation}
  \prodd(m_{i,j};x_i,x_j)\prodd(m_{i,j};x_j,x_i)^{-1}=1 \text{ for }i\ne
  j.
   \end{equation}
   (The right hand side of this equation will actually lie in the
   torus $T$, but we continue to designate the identity element by $1$
   because it is the identity of the nonabelian group $N$.)  We follow
   \cite[3.3]{rTits}. Let $x=\prodd(m_{i,j}-1;x_j,x_i)$, and let $k$
   denote $j$ or $i$, according to whether $m_{i,j}$ is even or odd.
   The desired relation can be rewritten
   \begin{equation} \label{FirstEquation}
     x x_j x^{-1} x_k^{-1}=1.
   \end{equation}
   Let $w$ be the image of $x$ in $W$. Then $w \srefl_j
   w^{-1}\srefl_k=1$ (see Section~\ref{GeneratorsAndRelations}), so that $xN_{\srefl_j}x^{-1}\subset
   N_{\srefl_k}$. Since $x_j\in N_{\srefl_j}$ and $x_k\in
   N_{\srefl_k}$ it follows that the left hand side of
   \eqref{FirstEquation} lies in $N_{\srefl_k}$.
   
   On the other hand, one can also rewrite  relation
   \eqref{FundamentalRelation} as
  \begin{equation}\label{SecondEquation}
  \prodd(m_{i,j};x_j,x_i)\prodd(m_{i,j};x_i,x_j)^{-1}=1 \text{ for }i\ne
  j.
  \end{equation}
  Let $l$ denote $i$ or $j$, according to whether $m_{i,j}$ is even or
  odd. The same reasoning as before shows that the
  left hand side of \eqref{SecondEquation} lies in $N_{\srefl_l}$. The
  conclusion is that the right hand side of \eqref{FirstEquation} lies
  in $N_{\srefl_k}\cap N_{\srefl_l}$. Since $G$ is simply connected, we
  are done by Proposition~\ref{IntersectionLemma}.
 \end{proof}

\section{Lattices, tori, and extensions at the prime~$2$}
\label{CMenagerieTwo}

In this section we copy the results of Sections~\ref{CMenagerie} and
\ref{CNormalizerExt}, replacing $\Z$ by the ring $\Ztwo$ of $2$--adic
integers; however, we do not formulate the notion of a root
system over $\Ztwo$. The proofs are omitted.

\begin{NumberedSubSection}{\Mrtlattices/}
  A \emph{\tlattice} \Lc\ is a fin\-ite\-ly generated free module over
  the ring $\Zt$ of $2$--adic integers.  A \emph{reflection} on
  $\Lc$ is an automorphism $\refl$ of \Lc/ which is
  conjugate in $\Aut(\Q\Tensor \Lc)\cong\GL(r,\Q_2)$ to a diagonal
  matrix $\diag(-1,1,\ldots,1)$.

\begin{defn}\label{tDefineMarking}
  Suppose that $\refl$ is a reflection on \atlattice/ $\Lc$. A
  \emph{strict marking} for $\refl$ is a pair $(b,\beta)$, where $b\in
  \Lc$ and $\beta\co \Lc\to\Ztwo$ is a homomorphism such that for any
  $x\in \Lc$, $ \refl(x) = x + \beta(x) b$.  Two strict markings
  $(b,\beta)$ and $(b',\beta')$ are \emph{equivalent} if $(b,\beta)=
  (ub,u^{-1}\beta)$ for a unit $u\in\Ztwo$.  A \emph{marking} for
  $\refl$ is an equivalence class $\{(b,\beta)\}$ of strict markings.
\end{defn}

\begin{rem}\label{tHowManyMarkings}
  As in Remark~\ref{HowManyMarkings}, a reflection has two markings if it is
  trivial mod~$2$, and one marking otherwise.
\end{rem}

Suppose that $\Lc$ is \atlattice, $\refl$ is a reflection on $\Lc$
with marking $\{(b,\beta)\}$, and $w$ is an automorphism of $\Lc$. Let
$w\cdot\{(b,\beta)\}$ denote the marking for $w\refl w^{-1}$ given by
$\{(w(b),\beta\circ w^{-1})\}$.

\begin{defn}\label{tDefineMarkedLattice}
  \emph{\Anrtlattice} is \atlattice/ $\Lc$ together with a finite
  subgroup $W$ of $\Aut(\Lc)$ which is generated
  by the reflections it contains. \emph{\Amrtlattice/} 
  is \anrtlattice/ $(\Lc,W)$ together with markings
  $\{(b_\refl,\beta_\refl)\}$, one for each reflection $\refl$ in
  $W$, such that for each element $w\in W$ and reflection
  $\sigma$,
  $
      w\cdot \{(b_\refl, \beta_\refl)\}= \{(b_{w\refl w^{-1}},
        \beta_{w\refl w^{-1}})\}
  $.
\end{defn}

\end{NumberedSubSection}

\begin{NumberedSubSection}{\Mrttori/}

Let $\dTorus(n)$ denote $(\Ztinfty)^n$. A \emph{\ttorus/} is a
discrete group isomorphic to $\dTorus(n)$ for some $n$. Any \ttorus/
\Td/ gives \atlattice/ $\Hom(\dTorus(1),\Td)$; conversely, \atlattice/
\Lc/ gives \attorus/ $\dTorus(1)\Tensor\Lc$. These two  constructions
are inverse to one another up to natural isomorphism, and induce an
equivalence between the category of \ttori/ and the category of \tlattices.

\begin{defn}\label{tDefineReflection}
  An automorphism $\refl$ of \attorus/ $\Td$ is a
  \emph{reflection} (respectively,  \emph{trivial} mod~$2$) if the induced
  automorphism of $\Hom(\dTorus(1),\Td)$ is a reflection (respectively,
  trivial mod~$2$).  An element $x$ of $\Td$ is said to be
  \emph{strongly $\refl$--negative} if $x$ lies in the maximal
  divisible subgroup $\Td_0^-(\refl)$ of $\Td^-(\refl)$.
\end{defn}

\begin{rem}\label{tDefineStronglyNegative}
  An automorphism $\refl$ of \Td/ as above is trivial mod~$2$ if and
  only if it acts as the identity on ${}_2T=\{x\in T\mid 2x=0\}$.
  If $T=\dTorus(1)\Tensor \Lc$, then $x\in \Td$ is strongly
  $\refl$--negative if and only if $x\in \dTorus(1)\Tensor \Lc^-(\refl)$.
\end{rem}

\begin{defn}\label{tDefineTorusMarking}
  Suppose that $\refl$ is  a reflection on \attorus/ \Td. A
  \emph{marking} for $\refl$ is an element $h\in\Td$ such that
  \begin{enumerate}
  \item $h$ is strongly $\refl$--negative (Remark~\ref{tDefineStronglyNegative}),
  \item $2h=0$, and
  \item $h\ne0$ if $\refl$ is nontrivial mod~$2$.
  \end{enumerate}
\end{defn}

\begin{defn}\label{tDefineMarkedReflTorus}
  \emph{\Anrttorus/}  is \attorus/ $\Td$ together with a finite
  subgroup $W$ of $\Aut(\Td)$ which is generated by the reflections it
  contains. \emph{\Amrttorus/} is \anrttorus/
  $(\Td,W)$ together with markings $h_\refl\in \Td$, one for each
  reflection $\refl$ in $W$, such that for  $w\in W$,
  $
          h_{w\refl w^{-1}} = w(h_\refl)
  $.
\end{defn}

\end{NumberedSubSection}

\begin{prop}\label{tInterpretMarkedTorus}
  If \Lc/ is \atlattice, there is a natural bijection between marked
  reflection structures on \Lc/ and marked reflection structures on
  the \ttorus/ $\dTorus(1)\Tensor\Lc$.
\end{prop}

\begin{NumberedSubSection}{The reflection extension}\label{tDiscussReflectionExt}
(Compare with Section~\ref{DReflectionExtension}.) Suppose that $(\Td,W)$ is \anrttorus, with associated
\rtlattice/ $\Lc$. Let $\refls$ be the
set of reflections in $W$, and write $\refls=\coprod_i\refls_i$ as a
union of conjugacy classes of reflections. For each index $i$ choose a
reflection $\crefl_i$ in $\refls_i$, let $C_i$ be the centralizer
of $\crefl_i$ in $W$, and let $a_i\in\Lc$ be a nonzero element with
$\crefl_i(a_i)=-a_i$. It is clear that each element of $C_i$ takes
$a_i$ to $\pm a_i$, and that $C_i$ is isomorphic to
$\Z/2 \times C_i\pperp$, where $C_i\pperp$ is the
subgroup of $C_i$ consisting of elements which fix $a_i$.
Consider the extension
\begin{equation}\label{tCanonicalExtension}
   1 \to \Z \to \Z\times C_i\pperp\to \Z/2\times
   C_i\pperp\to 1.
\end{equation}
By Shapiro's lemma, this gives rise to an extension $\RE_i(W)$ of $W$
by $\Z[\refls_i]$, and it is easy to see that up to canonical
isomorphism this extension of $W$ does not depend on the choice of
representative $\crefl_i$ for the conjugacy class $\refls_i$.

  \begin{defn}\label{tDefineReflectionExtension}
    The \emph{reflection extension} $\RE( W)$ of $W$ is the extension
    by $\Z[\refls]$ given by the sum of the extensions $\RE_i(W)$, in
    other words, by the fibre product of the extensions over~$W$.
  \end{defn}

\end{NumberedSubSection}

\begin{NumberedSubSection}{The normalizer extension}

Let $(\Td, W)$ be \anrttorus, and $\refls$ the set of
reflections in $W$. Given a marking $\{h_\refl\}$ for $(\Td, W)$,
ie a marked reflection torus $\mcT=(\Td,W,\{h_\refl\})$, there is a
$W$--map
$
    f^{\mathcal T}\co  \Z[\refls]\to \Td
$
which sends a reflection $\refl$ to $h_\refl$.

  \begin{defn}\label{tDefineNormalizerExt}
    Suppose that $\mcT=(\Td,W,\{h_\refl\})$ is \amrttorus.  The
    \emph{normalizer extension} $\Normd(T,W,\{h_\refl\})$ of $W$ by
    $\Td$ is the image under $f^{\mcT}\co \Z[\refls]\to \Td$ of the
    reflection extension $\RE(W)$.
\end{defn}

\end{NumberedSubSection}

\section{Classification of marked complete reflection lattices}\label{CtLattices}

In this section we prove Proposition~\ref{ClassifyCompleteLattices}, which roughly
states that, one example aside, every \mrtlattice/ arises from a
marked reflection lattice, ie from a classical root system.

It will be convenient to consider pairs $(A,W)$ in which $A$ is an
abelian group and $W$ is a group acting on $A$; one example is a
reflection lattice $(L,W)$. It is clear what it means for two such
objects to be isomorphic. The product $(A,W)\times (A', W')$ of two is
the pair $(A\times A', W\times W')$, with $W\times W'$ acting on
$A\times A'$ in the obvious way.  Note that the product of two
complete and/or marked reflection lattices is naturally a reflection
lattice of the same type.  If $(L,W)$ is a marked reflection lattice,
then $(\Ztwo\Tensor L, W)$ is naturally \amrtlattice.

\begin{defn}\label{DefineCoxeterType}
  A marked complete reflection lattice is \emph{of Coxeter type} if
  it is isomorphic to $(\Ztwo\Tensor L, W)$ for some marked
  reflection lattice $(L,W)$.
\end{defn}

Let $(\LD,\WD)$ be the complete reflection lattice of rank~$3$
associated to the exceptional $2$--compact group $\DI(4)$
\cite{rDWnew,rNotbohmDIFour}.

\begin{rem}\label{DIFourStuff}
  The group $\WD$ is isomorphic to $\Z/2\times \GL(3,\Ftwo)$; the
  central $\Z/2$ in \WD/ acts by $-1$ on $\LD$, and \WD/ acts on
  $\Z/2\Tensor\LD\cong(\Z/2)^3$ via the natural representation of
  $\GL(3,\Ftwo)$.  These last two facts essentially determine the
  action of $\WD$ on $\LD$ \cite[Section~4]{rDWnew}. Every reflection in $\WD$ is nontrivial mod~$2$,
  and so has a unique marking (Remark~\ref{tHowManyMarkings}); in particular,
  there is a unique way to give $(\LD, \WD)$ the structure of
  \mrtlattice.
\end{rem}

\begin{defn}\label{DefineTypeDIfour}
  A marked \rtlattice/ is \emph{of type $\DI(4)$}
  if it is isomorphic to a product of copies of $(\LD,\WD)$.
\end{defn}

\begin{prop}\label{ClassifyCompleteLattices}
  Every \mrtlattice/  $(\Lc,W)$ is isomorphic to
  a product $(\Lc_1,W_1)\times(\Lc_2,W_2)$, in which the first factor
  is of Coxeter type and the second is of type $\DI(4)$.
\end{prop}

This depends on a few lemmas.

\begin{defn}
  \Amrtlattice/ $(\Lc,W)$ is \emph{rationally of Coxeter type} if
  there is a reflection lattice $(L,W)$ such that
  $(\Q\Tensor\Lc, W)$ is isomorphic as a reflection lattice to $(\Qtwo\Tensor L, W)$.
\end{defn}

\begin{lem}\label{Decomplete}
  Suppose that $(\Lc,W)$ is \amrtlattice/  which is
  rationally of Coxeter type. Then $(\Lc,W)$ is of Coxeter
  type.
\end{lem}

\begin{proof}
  Since $(\Lc,W)$ is rationally of Coxeter type, there is a lattice
  $L'$ and an action of $W$ on $L'$ such that $(\Q\Tensor\Lc, W)$
  is isomorphic to $(\Qtwo\Tensor L', W)$; we use this isomorphism to
  identify $\Q\Tensor\Lc$ with $\Qtwo\Tensor L'$.  
  Under this isomorphism,  $\Lc'=\Ztwo\Tensor L'$ is \atlattice/
  in $\Q\Tensor\Lc$ preserved by $W$. Adjust $L'$ by a
  power of $2$ so that $\Lc\subset\Lc'$. Now define $L$ by the
  pullback diagram:
  \[
  \begin{CD}
      L @>>> \Lc \\
      @VVV   @VVV \\
      L' @>>> \Lc'
  \end{CD}
  \]
  The group $\Lc'/\Lc$ is a finite $2$--group, and the map
  $L'\to\Lc'\to\Lc'/\Lc$ is surjective, this last by Nakayama's lemma
  and the fact that $\Z/2\Tensor L'\iso\Z/2\Tensor\Lc'$. It follows
  that $L'/L$ is isomorphic to $\Lc'/\Lc$, and that $L$ is a
  lattice in $\Q\Tensor L'$ with $\Ztwo\Tensor L\iso\Lc$.  In
  particular, $(\Lc,W)$ is isomorphic to $(\Ztwo\Tensor L,W)$.
  The fact that under this isomorphism the  marking on $\Lc$ is
  derived from a marking on $L$  follows from the observation that
  possible markings of $(\Lc,W)$ correspond bijectively to markings of
  $(L,W)$.
\end{proof}

\begin{lem}\label{LSplitOffDIfour}
  Suppose that $(\Lc,W)$ is \anrtlattice, such that
  $(\Q\Tensor\Lc, W)$ decomposes as a product $(V_1,W_1)\times(\Q\Tensor
  \LD,\WD)$. Then $(\Lc,W)$ decomposes as the product $(\Lc\cap
  V_1,W_1)\times(\LD,\WD)$ of \rtlattices. If
  $(\Lc,W)$ is marked, then there is a unique way to mark the two
  factors of the product in such a way that the product decomposition
  respects the markings.
\end{lem}

\begin{proof}
  Write $\Q\Tensor\Lc\iso V_1\times V_2$ and $W\iso W_1\times W_2$,
  where $V_2\iso\Q\Tensor\LD$ and $W_2\iso\WD$. Let $\Lc_i=\Lc\cap
  V_i$, so that $\Lc_i$ is \atlattice/ which is a module over $W_i$
  and $\Q\Tensor \Lc_i\iso V_i$. By elementary rank considerations, if
  $\refl\in W$ is a reflection, then either $\refl\in W_1$ or
  $\refl\in W_2$. Consequently, $W_1$ and $W_2$ are generated by the
  reflections they contain, and both $(\Lc_1, W_1)$ and $(\Lc_2,W_2)$
  are \rtlattices. Now adjust the copy of $\LD$ in
  $V_2$ by multiplication by a power of $2$ so that $\LD\subset \Lc_2$
  but $\LD$ is not contained in $2\Lc_2$. The induced map $f\co 
  \Z/2\Tensor\LD\to\Z/2\Tensor\Lc_2$ is then a nontrivial map between
  $W_2$--modules which are $\Ftwo$--vector spaces of the same dimension;
  since $\Z/2\Tensor\LD$ is simple as a $W_2$--module
  (Remark~\ref{DIFourStuff}), the map $f$ is an isomorphism and it follows
  from Nakayama's lemma that $\Lc_2$ is isomorphic to $\LD$ as a
  $W_2$--module. Let $\bar L=\Z/2\Tensor\Lc$ and $\bar L_i=\Z/2\Tensor
  \Lc_i$. Consider the two $W$--module maps
  \[
     g\co  \Lc_1\times\Lc_2\to \Lc \quad\quad
     g'\co \bar L_1\times\bar L_2\to\bar L.
  \]
  Again by Nakayama's lemma, in order to prove that $g$ is an
  isomorphism it is enough to prove that $g'$ is an isomorphism, or
  even, by a rank calculation, that $g'$ is injective. Let $g'_i$
  denote the natural map $\bar L_i\to\bar L$ and let $K$ be the kernel
  of $g'$. It is clear from the definitions that $g'_i$ is injective
  (ie, if $x\in \Lc$ and $2x\in\Lc_i$ then $x\in\Lc_i$), so that the
  two maps $K\to\bar L_1$ and $K\to\bar L_2$ are injective. Since
  $W_2$ acts trivially on $\bar L_1$, $W_2$ must act trivially on $K$,
  and so the fact that $K=\{0\}$ follows from the fact that $\bar
  L_2\iso\Z/2\Tensor\LD$ has no submodules with a trivial
  $W_2$--action. The quickest way to obtain the statement about
  markings is to interpret the markings as elements of the associated
  discrete tori (\ref{tInterpretMarkedTorus}).
\end{proof}

\begin{prop}\label{ListTheReflectionGroups}
  The only irreducible finite reflection groups over $\Q_2$ are the
  classical Weyl groups, with their standard reflection representations, and
  the group $\WD$, with the reflection representation derived from $\DI(4)$. 
\end{prop}

\begin{proof}
  This follows from \cite{rClarkEwing}. Clark and Ewing observe that
  the Schur index of a reflection representation is~$1$, so that the
  representation is defined over its character field, that is, the field
  extension of $\Q$ generated by the characters of the group elements.
  They then determine the character fields of all of the irreducible
  complex reflection groups. The only character fields contained in
  $\Q_2$ are the fields derived from classical Weyl group
  representations (which are all $\Q$) and the field derived from the
  reflection representation of $\WD$ (which is $\Q(\sqrt{-7})$).  The
  complex reflection representation of $\WD$ is unique up to
  conjugacy, which easily implies that the associated 
  $\Q_2$--reflection representation is also unique.
\end{proof}

\begin{proof}[Proof of Proposition~\ref{ClassifyCompleteLattices}]
  Let $(\Lc,W)$ be \anrtlattice.  According to the classification of
  $2$--adic rational reflection groups (see Proposition~\ref{ListTheReflectionGroups}), we can
  write $(\Q\Tensor\Lc,W)$ as a product $(V_1, W_1)\times(
  \Q\Tensor\LD ,\WD)^k$ for some $k\ge0$, where $W_1$ is a classical
  Weyl  group. Repeated applications of Lemma~\ref{LSplitOffDIfour} give
  that $(\Lc,W)$ is equivalent as \amrtlattice/ to
  $(\Lc_1,W_1)\times(\LD,\WD)^k$, where $(\Lc_1,W_1)$ is rationally of
  Coxeter type. The result follows from Lemma~\ref{Decomplete}.
\end{proof}

\section{Completions of compact Lie groups}\label{CPassingToTwo}

In this section we give some results about completions of Lie groups
that are used in Section~\ref{CTCGroups}.  Suppose that $G$ is a connected
compact Lie group of rank~$r$ and that $X$ is the $2$--compact group
obtained by taking the $2$--completion of $G$ (so that $BX$ is the
$2$--completion of $BG$).  Let $(T,W,\{h_\refl\})$ be the marked
reflection torus obtained from $G$, and $\Td$ the $2$--primary torsion
subgroup of $T$ \eqref{TwoPrimary}, so that $\Td$ is \attorus/ of
rank~$r$. It is easy to see that $(\Td, W,\{h_\refl\})$ is \amrttorus.
Let $\Norm(W)$ be the normalizer extension of $W$ by $T$ determined by
$(T,W,\{h_\refl\})$ (Definition~\ref{DefineNormalizerExt}) and $\Normd(W)$ the
parallel normalizer extension of $W$ by $\Td$
(Definition~\ref{tDefineNormalizerExt}). Essentially by construction there is a
homomorphism $\Normd(W)\to\Norm(W)$ which lies in a map of exact
sequences:
  \[
  \begin{CD}
     1 @>>> \Td @>>> \Normd(W) @>>> W @>>> 1 \\
     @.  @V\subset VV  @VVV @V\|VV @.\\
     1 @>>> T @>>> \Norm(W) @>>> W @>>> 1
  \end{CD}
  \]
  By Theorem~\ref{IdentifyTorusNormalizer}, the inclusion $T\to G$ extends to
  a homomorphism $\Norm(W)\to G$.  Combining these homomorphisms with
  the completion map $BG\to BX$ gives composite maps $B\Td\to BT\to
  BG\to BX$ and $B\Normd(W)\to B\Norm(W)\to BG\to BX$, which, in the
  language of \cite{rDWmethods}, correspond to homomorphisms $\Td\to
  X$ and $\Normd(W)\to X$.

  \begin{lem}\label{IdentifyTandNT}
    In the above situation, $\Td\to X$ is a  maximal \ttorus/ for
    $X$, $\Normd(W)$ is the normalizer of $\Td$ in $X$, and
    $\Normd(W)\to X$ is the natural homomorphism.
  \end{lem}

  See Aside~\ref{DiscussDiscreteTorus} for a description of the maximal
  \ttorus/ in a $2$--compact group.

  \begin{proof}[Proof of Lemma~\ref{IdentifyTandNT}]
    Consider first the chain of maps
    \[
        B\Td \to BT \to BG \to BX.
    \]
    The left hand map induces an isomorphism on mod~$2$ homology (see
    \cite[Section~6]{rDWmethods}), as does the completion map $BG\to BX$.
    Let $X/\Td$ denote the homotopy fibre of $B\Td\to BX$, $G/T$ the
    homotopy fibre of $BG\to BX$, and $G/\Td$ the homotopy fibre of
    $B\Td\to BG$. Let $\hat T$ denote the closure \cite[Section~3]{rDWcenter} of
    $\Td$ (ie $B\Tc$ is the $2$--completion of $B\Td$), so that
    $B\Td\to BX$ extends canonically to $B\hat T\to BX$, and denote
    the homotopy fibre of this last map by $X/\hat T$. It follows that
    all of the maps in the diagram
    \[
       X/\hat T\leftarrow    X/\Td\leftarrow G/\Td \to G/T
    \]
    induce isomorphisms on mod~$2$ homology. In particular, $X/\hat T$
    has finite mod~$2$ homology, and so $\hat T\to X$ is a
    monomorphism \cite[3.2, 9.11]{rDWmethods}. It is straightforward
    to argue that the Euler characteristic $\chi_{\Qtwo}(X/\hat T)$
    \cite[1.5]{rDWmethods} is the same as the rational Euler
    characteristic $\chi_{\Q}(G/T)$; since this last is nonzero
    \cite[5.3]{rDWelementary}, the map $\hat T\to X$ is a maximal
    torus for $X$ \cite[2.15]{rDWcenter}. Since $\Td\to \hat T$ is a
    discrete approximation, $\Td\to X$ is a \td/ maximal torus for
    $X$. The statements about $\Normd(W)$ amount to an identification
    of $W$ with the Weyl group $W_X$ of $X$. There is certainly a map
    $W\to W_X$, since $W$ maps to the self-equivalences of $BT$ over
    $BG$ and hence by completion to the self-equivalences of $B\hat T$
    over $BX$ (see \cite[9.5]{rDWmethods}). The map is a monomorphism
    because the completion map $\GL(r,\Z)\iso
    \Aut(\pi_1T)\to\Aut(\pi_1\hat T)\iso\GL(r,\Ztwo)$ is injective,
    and consequently an isomorphism by counting:
    $\#(W)=\chi_{\Q}(G/T)$ \cite[5.3]{rDWelementary}, while
    $\#(W_X)=\chi_{\Qtwo}(X/\hat T)$ \cite[9.5]{rDWmethods}.
  \end{proof}

We will also need to refer to the following results.

\begin{prop}\label{SameConjugacies}
  Suppose that $G$ is a connected compact Lie group, that $X$ is
  the $2$--completion of $G$, and that $\pi$ is a finite $2$--group.
  Then the map which sends a homomorphism $f\co \pi\to G$ to the composite
  $B\pi\RightArrow{Bf} BG\to BX$ induces a bijection between conjugacy
  classes of homomorphisms $\pi\to G$ and homotopy classes of maps
  $B\pi\to BX$. 
\end{prop}
 
\begin{proof}
  This amounts to combining \cite[1.1]{rDZ} with the fact that  the
  completion map $BG\to BX$ induces a bijection from homotopy classes
  of maps $B\pi\to BG$ to homotopy classes of maps $B\pi\to BX$. This
  last follows from \cite[2.4]{rDZ}.
\end{proof}

\begin{prop}\label{SameCentralizers}
  Suppose that $G$ is a connected compact Lie group with maximal torus
  $T$, and that $\Sd\subset T$ is a discrete torus with topological
  closure $\Scl$. Let $X$ be the $2$--completion of $G$, $G_{\Sd}$ the
  centralizer of $\Sd$ in $G$, $G_{\Scl}$ the centralizer of $\Scl$ in
  $G$, and $X_{\Sd}$ the centralizer of $\Sd$ in $X$. Then the natural
  map $G_{\Scl}\to G_{\Sd}$ is an isomorphism, while the natural map
  $BG_{\Sd}\to BX_{\Sd}$ is a $2$--completion map.
\end{prop}

\begin{proof}
  It is elementary that $G_{\Scl}\to G_{\Sd}$ is an isomorphism. The
  fact that $BG_{\Sd}\to BX_{\Sd}$ is a $2$--completion map follows
  from the main theorem of \cite{rNotbohmMaps} and the argument in the
  proof of \cite[2.4]{rDZ}. Note that $X_{\Sd}$ is $2$--complete by
  \cite[5.7, 6.21]{rDWmethods}.
\end{proof}
\section{$2$--compact groups}\label{CTCGroups}

In this section we associate \amrttorus/ to a connected $2$--compact
group $X$ (Proposition~\ref{FundamentalConjugacies}), and show that the associated
normalizer extension of the Weyl group is isomorphic to the normalizer
of a maximal \ttorus/ in $X$ (Proposition~\ref{IdentifyNT}).

\begin{aside}\label{DiscussDiscreteTorus}
  A \emph{complete torus} $\hat T$ is a $2$--compact group with the
  property that its classifying space $B\hat T$ is equivalent to the
  $2$--completion of the classifying space of an ordinary torus. Any
  $2$--compact group $X$ has a \emph{maximal (complete) torus}, which
  is technically a complete torus $\hat T$ together with a
  monomorphism $\hat T\to X$ with the property that the Euler characteristic of
  $X/\hat T$ is nonzero \cite[Section~8]{rDWmethods}. For such a maximal torus, the space
  $\mathcal{W}$ of self-maps of $B\hat T$ over $BX$ is homotopically
  discrete, and its monoid of components is a finite group $W$ called
  the \emph{Weyl group} of $X$. The monoid $\mathcal W$ acts on $B\Tc
  $, and the Borel construction of this action is denoted $BN(\hat T)$; this lies in a fibration sequence
  \[  B\hat T\to BN(\hat T)\to BW.\]
   The loop space of $BN(T)$ is called the \emph{normalizer of  $\hat T$
    in $X$}.

  Any complete torus has a \emph{discrete approximation}, which is a
  discrete torus $\Td$ together with a homomorphism $\Td\to\hat T$
  with the property that the induced map $B\Td\to B\Tc$ induces an
  isomorphism on mod~$2$ homology. If $\Tc\to X$ is a maximal torus
  for $X$, then the composite $\Td\to\Tc\to X$ is called a
  \emph{maximal discrete torus} for $X$.  See \cite[3.13]{rDWcenter}
  for an explanation of how in this case the approximation $\Td\to\Tc$ can be
  extended to give a discrete approximation  $\Nd$ for $N(\Tc)$: $\Nd$ is a
  discrete group which lies in a short exact sequence
  \[
      \{1\} \to \Td\to\Nd\to W\to\{1\},
  \]
  and it is called the \emph{normalizer of $\Td$} in $X$.
\end{aside}

Suppose that $X$ is a connected $2$--compact group of rank~$r$ with
maximal \ttorus/ \Td, Weyl group $W$, and discrete torus normalizer
\Nd. Then $\Td\iso \dTorus(r)$ and the conjugation action of $W$ on
$\Td$ identifies $W$ with a subgroup of $\Aut(\Td)$ generated by
reflections \cite[9.7]{rDWmethods}. In particular, $(\Td,W)$ is a \td/
reflection torus.

For each reflection $\refl\in W$ we will construct a subgroup
$\Nd_\refl$ of $\Nd $, called the \emph{root subgroup} associated to
$\refl$. We will also single out a specific element $h_\refl\in\Td$,
the \emph{marking} associated to $\refl$.  The technique is to mimic
in homotopy theory the constructions of Section~\ref{CCompactLie}.

Given $\refl$, let $\Td^+_0(\refl)$ denote the maximal divisible
subgroup of $\Td$ on which $\refl$ acts trivially. It is easy to see
that $\Td^+_0(\refl)$ is \attorus/ of rank~$r-1$; if $\Lc$ is
the \tlattice/  corresponding to \Td, then under the isomorphism
$\Td\iso\dTorus(1)\Tensor\Lc $, $\Td^+_0(\refl)$ is given by
$\dTorus(1)\Tensor\Lc^+(\refl)$. Let $\RX\refl$ denote the centralizer
of $\Td^+_0(\refl)$ in $X$.

\begin{aside}\label{MappingGivesCentralizer}
  The inclusion homomorphism $i\co \Td^+_0(\refl)\to\Td\to X$ is by
  definition represented by a map $Bi\co B\Td^+_0(\refl)\to BX$ of
  classifying spaces. Note that the domain of $Bi$ is a space of type
  $K(\pi,1)$. The centralizer $\RX\refl$ is by definition the
  $2$--compact group whose classifying space $B\RX\refl$ is the mapping
  space component $\Map(B\Td^+_0(\refl),BX)_{Bi}$.
\end{aside}

Then $\RX\refl$ is a connected $2$--compact subgroup of $X$
\cite[7.8]{rDWcenter} with maximal discrete torus $\Td$, and Weyl group given
by the subgroup $\{1,\refl\}$ of $W$ \cite[7.6]{rDWcenter}. The
normalizer of $\Td$ in $\RX\refl$ can be identified with the algebraic
centralizer of $\Td^+_0(\refl)$ in $\Nd$; we denote this group by
$\Nd(\refl)$. Note that $\Nd(\refl)$ is an extension of $\Z/2$ by
\attorus, and so is a $2$--discrete toral group \cite[6.5]{rDWmethods}.

\begin{aside}
  Let $A=B\Td^+_0(\refl)$.
  The classifying space maps
  \[
  A=B\Td^+_0(\refl)\RightArrow{u} B\Td \RightArrow{v} B\Nd \RightArrow{w} BX
  \]
induce maps
\begin{equation}\label{SourceOfRootSubgroups}
  \Map(A,B\Td)_{u}\to \Map(A,B\Nd)_{vu}\to\Map(A,BX)_{wvu},
\end{equation}
The mapping space component on the far right is $B\RX\refl$. By
covering space theory, $\Map(A,B\Td)_u$ is canonically equivalent to
$B\Td$, the equivalence being given by evaluation at the basepoint of
$A$. Similarly, $\Map(A,B\Nd)_{vu}$ is the classifying space of the
algebraic centralizer of $\Td_0^+(\refl)$ in $\Nd$.  It follows from
\cite[7.6]{rDWcenter} that the composite map in
\eqref{SourceOfRootSubgroups} represents the inclusion in $\RX\refl$ of
a maximal \ttorus, and that the right-hand map represents the
inclusion of the normalizer of this torus.
\end{aside}

The only simply-connected $2$--compact group of rank~$1$ is the
two-completion \SUc/ of $\SU(2)$ (see  \cite{rDMWuniqueness}); the
center of $\SUc$ is $\Z/2$ and the adjoint form is the $2$--completion
$\SOc$ of $\SO(3)$.  Given this, it follows from the product splitting
theorem for $p$--compact groups \cite{rDWproduct} that $\RX\refl$ is
given by one of the following three possibilities (see
\cite[7.7]{rDWcenter}, \cite{rNMcenter}), where $\cTorus(k)$ denotes a
$2$--complete torus of rank~$k$ and $\Uc$ is the $2$--completion of
$\U(2)$:
\begin{enumerate}
\item  A product $\SUc\times\cTorus(r-1)$
\item A product $\SOc\times\cTorus(r-1)$
\item A product  $\Uc\times\cTorus(r-2)$
\end{enumerate}   

Recall that $\Td^-_0(\refl)$ denotes the maximal divisible subgroup of
$\Td^-(\refl)$; this group is \attorus/ of rank~$1$.

\begin{defn}\label{tDefineRootSubgroup}
  The \emph{root subgroup} $\Nd_\refl$ is given by
\[
  \Nd_\refl=\{x\in \Nd(\refl)\mid x\text{ is conjugate in
    }\RX\refl\text{ to some }y\in \Td^-_0(\refl)\}.
\]
\end{defn}

\begin{aside}\label{HomotopicalConjugacy}
  This definition is couched in the language of \cite{rDWmethods}, but
  what it means concretely is the following.  Suppose that
  $x\in\Nd(\refl)$ has order $2^n$; note that by
  \cite[6.19]{rDWmethods}, $\Nd(\refl)$ is a union of finite
  $2$--groups.  Represent $x$ by $\chi(x)\co \Z/2^n\to \Nd(\refl)$. Then
  $x\in\Nd_\refl$ if and only if the composite
  \[
     B\Z/2^n\RightArrow{B\chi(x)} B\Nd(\refl)\to B\RX\refl
  \]
  lifts up to homotopy to a map $B\Z/2^n\to B\Td^-_0(\refl)$.
\end{aside}

Given that $\RX\refl$ is the $2$--completion of a compact Lie group,
the following statement comes from combining
Lemma~\ref{RootSubgroupIsSubgroup} with the discussion in
Section~\ref{CPassingToTwo}. Lemma~\ref{LXMarkings} is proved similarly. 

\begin{lem}\label{tRootSubgroupIsSubgroup}
  The set $\Nd_\refl$ is a subgroup of $\Nd(\refl)$ isomorphic to an
  extension of $\Z/2$ by a rank~$1$ \ttorus. More precisely,
  $\Nd_\refl\cap\Td=\Td^-_0(\refl)$, and the natural map
  $\Nd_\refl/\Td^-_0(\refl)\to\Nd(\refl)/\Td\iso\Z/2$ is an
  isomorphism. If $x, y\in \Nd_\refl\setminus \Td^-_0(\refl)$, then
  $x^2=y^2 \in \Td^-_0(\refl)$.
\end{lem}

\begin{defn}\label{DefineLXMarkings}
  The \emph{marking} $h_\refl\in \Td^-_0(\refl)$ is defined to be the
  image of $\Nd_\refl\setminus \Td^-_0(\refl)$ under the squaring map
  $x\mapsto x^2$.
\end{defn}

\begin{lem}\label{LXMarkings}
  The element $h_\refl\in \Td$ is a marking for $\refl$ in the sense of
  Definition~\ref{tDefineTorusMarking}.
\end{lem}

\begin{prop}\label{FundamentalConjugacies}
  Suppose that $X$ is a connected $2$--compact group with maximal
  \ttorus/ $\Td$, discrete torus normalizer \Nd, Weyl group $W$, root
  subgroups $\{\Nd_\refl\}$, and markings $\{h_\refl\}$ as above.
  Suppose that $x\in\Nd$ has image $w\in W$.  Then for any reflection
  $\refl\in W$,
  \[
  \begin{aligned}
  x\Nd_{\refl\vphantom{w^{-1}}} x^{-1} &= \Nd_{w\refl w^{-1}}\\
  x h_{\refl\vphantom{w^{-1}}}x^{-1} &= h_{w\sigma w^{-1}}
  \end{aligned}
  \]
  In particular, $(\Td,W,\{h_\refl\})$ is \amrttorus.
\end{prop}

\begin{proof}
  As with Proposition~\ref{LGConjugationInvariant}, the idea here is that
  conjugation with $x$ gives a symmetry of the whole situation, but
  there is something to check. Let $c$ denote conjugation with $x$,
  and let $\tau=w\refl w^{-1}$. All of the arrows in the following
  commutative diagram are clear, except for the right hand vertical
  one:
  \[
   \begin{CD}
    B T^+_0(\refl)@>>> B\Nd @>>> X \\
   @V Bc VV                   @V Bc VV  @V Bc' VV \\
    B T^+_0(\tau) @>>> B\Nd @>>> X
   \end{CD}
  \]
  If the right hand arrow exists, then taking mapping spaces (Aside~\ref{MappingGivesCentralizer}) gives a
  commutative diagram
  \[
   \begin{CD}
     B\Nd(\refl) @>>> B\RX\refl\\
       @V Bc V V        @VVV \\
     B\Nd(\tau) @>>> B\RX\tau
   \end{CD}
  \]
  which, in view of Definitions~\ref{tDefineRootSubgroup} and
  \ref{DefineLXMarkings}, yields the desired result. The issue then is
  to extend the self equivalence $Bc$ of $B\Nd$ to a self equivalence $Bc'$ of $BX$.
  One way to achieve this is to represent the map $\Nd\to X$ by a map
  $U\to V$ of simplicial groups; the map $c$ can then be realized by
  conjugation with a vertex of $U$, and it extends to an automorphism
  of $V$ given by conjugation with the image vertex. Another approach
  is to treat $c$ as part of the conjugation action of $\Nd$ on
  itself. Taking classifying spaces translates this into a
  basepoint-preserving action of $\Nd$ on $B\Nd$, an action which is
  homotopically captured by the associated sectioned fibration over
  $B\Nd$ with fibre $B\Nd$.  It is easy to see that this fibration is
  just the product projection $B\Nd\times B\Nd\to B\Nd$, with section
  given by the diagonal map.  This fibration evidently extends to the
  product fibration $BX\times BX\to BX$.
\end{proof}

\begin{rem}\label{NormalizerConjecture}
  It is natural to conjecture that the above \mrttorus, or
  equivalently the associated \mrtlattice/
  (Proposition~\ref{tInterpretMarkedTorus}), determines $X$ up to equivalence.
\end{rem}

The main result of this section is the following one.

\begin{prop}\label{IdentifyNT}
  Suppose that $X$ is a connected $2$--compact group, and that
  $(\Td,W,\{h_\sigma\})$ is the \mrttorus\ associated to $X$
  (see Proposition~\ref{FundamentalConjugacies}).  Then the normalizer $\Nd$ of $\Td$
  in $X$ is isomorphic to the normalizer extension $\Normd(\Td, W,\{h_\refl\})$
  (Definition~\ref{tDefineNormalizerExt}).
\end{prop}

The proof breaks up into two cases.  A connected $2$--compact
group $X$ is said to be \emph{of Coxeter type} (respectively, \emph{of
type $\DI(4)$}) if the associated complete marked reflection lattice
(Propositions~\ref{FundamentalConjugacies} and \ref{tInterpretMarkedTorus})
is of Coxeter type (respectively, of type $\DI(4)$); see
Definitions~\ref{DefineCoxeterType} and \ref{DefineTypeDIfour}. Note
that this is a property of the underlying \rtlattices, or equivalently
of the action of $W$ on $\Td$; the markings come along automatically
(see Remark~\ref{DIFourStuff} and the proof of Lemma~\ref{Decomplete}).

\begin{prop}\label{SplitIntoCases}
  Any connected $2$--compact group $X$ can be written as a product $X_1\times
  X_2$ of $2$--compact groups, where $X_1$ is of Coxeter type and $X_2$
  is of type $\DI(4)$.
\end{prop}

\begin{proof}
  This is a consequence of Proposition~\ref{ClassifyCompleteLattices} and
  \cite[6.1]{rDWproduct}. 
\end{proof}

A connected $2$--compact group $X$ is said to be \emph{semisimple} if
its center is finite (equivalently, if $\pi_1X$ is finite). If $X$ is
any connected $2$--compact group, the universal cover $\tilde X$ is
semisimple, and it follows from \cite{rDWproduct} and \cite{rNMcenter}
that $X$ is obtained from the product of $\tilde X$ and a $2$--complete
torus by dividing out \cite[8.3]{rDWmethods} by a finite central
subgroup.  We express this by saying that $X$ is a \emph{central
  product} of $\tilde X$ and a torus. The \emph{semisimple rank} of
$X$ is the rank of $\tilde X$; or equivalently the rank of $X$ minus
the rank of the center of $X$. There are analogous notions for compact
Lie groups.

The rank~1 case of the following statement follows from
\cite{rDMWuniqueness}; the rank~2 case is proved by Bauer et al
cite[6.1]{rKN}.

\begin{prop}\label{RankTwo}
  If $X$ is a semisimple $2$--compact group of rank $\le 2$,
  then $ X$ is the $2$--completion of a connected compact
  Lie group.
\end{prop}

If $G$ is a semisimple Lie group, the center of $G$ is a finite
abelian group $\mathcal{Z}(G)$. It is easy to see that the center of
the corresponding $2$--compact group $\hat G$ is the quotient of
$\mathcal{Z}(G)$ by the subgroup of elements of odd order; for
instance, combine the parallel calculations \cite[7.5]{rDWcenter} and
\cite[8.2]{rDWelementary} with Lemma~\ref{IdentifyTandNT}. This implies that
any central product of $\hat G$ with a $2$--complete torus is the
$2$--completion of some central product of $G$ with a torus. In
conjunction with Proposition~\ref{RankTwo}, this observation gives the following.

\begin{prop}\label{SemisimpleRankTwo}
  If $X$ is a connected $2$--compact group of semisimple rank at most 2,
  then $X$ is the $2$--completion of a connected compact Lie group.
\end{prop}

Suppose that $X$ is a connected $2$--compact group with \mrttorus/
$(\Td,W,\{h_\refl\})$ and discrete torus normalizer $\Nd$ as above.
Let $A$ be a subgroup of $\Td$, $X_A$ the centralizer of $A$ in $X$,
and $W_A$ the subgroup of $W$ consisting of all elements which leave
$A$ pointwise fixed. Then $\Td$ is a maximal \ttorus/ for $X_A$, and
$W_A$ is the Weyl group of $X_A$. The normalizer $\Nd_A$ of $\Td$ in
$X_A$ is the centralizer of $A$ in $\Nd$, or equivalently the pullback
over $W_A\subset W$ of the surjection $\Nd\to W$. We assume that $X_A$
is connected; this is always true if $A$ is \attorus, and the general
criterion is in \cite[7.6]{rDWcenter}.  Under this assumption, $W_A$
is generated by the reflections it contains \cite[7.6]{rDWcenter}, and
so $(\Td,W_A)$ is \anrttorus.  Let $\{h_\refl\}_{\refl\in W_A}$
denote the set of markings for the reflections in $W_A$ obtained from
the given set $\{h_\refl\}$ of markings for the reflections in $W$.
It is clear that $(\Td, W_A, \{h_\refl\}_{\refl\in W_A})$ is a marked
reflection torus.  For each reflection $\sigma\in W_A$, let
$\Nd_\refl$ be the root subgroup for $\refl$ computed in $X$, and
$\Nd_{A,\refl}\subset \Nd_A$ the root subgroup for $\refl$ computed in
$X_A$. Given the inclusion $\Nd_A\subset\Nd$, both of these root
subgroups can be considered as subgroups of $\Nd$.

\begin{lem}\label{CutDownToCentralizer}
  In the situation described above, $(\Td, W_A, \{h_\refl\}_{\refl\in W_A})$ is
  the discrete marked reflection torus associated to $X_A$ by
  Proposition~\ref{FundamentalConjugacies}. Moreover, for each reflection $\refl\in W_A$,
  $\Nd_{A,\refl}=\Nd_\refl$ as a subgroup of $\Nd$.
\end{lem}

\begin{proof}
  By Definition~\ref{DefineLXMarkings}, the first statement follows from the
  second. Let $B=\Td^+_0(\refl)$, $\RX\refl$ the centralizer of $B$ in
  $X$, and $\RXA\refl$ the centralizer of $B$ in $X_A$. Both
  $\RX\refl$ and $\RXA\refl$ are connected \cite[7.8]{rDWcenter}, and
  the monomorphism \cite[5.2, 6.1]{rDWmethods} $X_A\to X$ induces a
  monomorphism $\RXA\refl\to\RX\refl$. The monomorphism is an
  equivalence by \cite[4.7]{rDWcenter}, since both $\RXA\refl$ and
  $\RX\refl$ have Weyl group $\Z/2$ generated by $\refl$. A similar
  argument shows that the centralizer $\Nd_A(\refl)$ of $B$ in $\Nd_A$ is
  the same as the centralizer $\Nd(\refl)$ of $B$ in $\Nd$: both groups
  have $\Td$ as a maximal normal divisible abelian subgroup, and the quotient of
  each group by $\Td$ is $\{1,\refl\}$. The second
  statement of the lemma thus follows from Definition~\ref{tDefineRootSubgroup}
\end{proof}

\begin{lem}\label{SameRootSubgroups}
  Suppose that $G$ is a connected compact Lie group, and $X$ the
  $2$--compact group obtained as the $2$--completion of $G$. Let
  $W$ be the Weyl group of $G$ (or of $X$), $N$ the normalizer of
  a maximal torus in $G$, $\Nd$ the normalizer of a maximal
  discrete torus $\Td$ in $X$, and $N_\refl$ (respectively, $\Nd_\refl$) the
  root subgroup of $N$ (respectively, $\Nd$) corresponding to a reflection
  $\refl\in W$.  Then under the identification given by
  Lemma~\ref{IdentifyTandNT} of $\Td$
  as a subgroup of $T$ and $\Nd$ as a subgroup of $N$,
  $\Nd_\refl=\Nd\cap N_\refl$.
\end{lem}

\begin{proof}
  Let $T^+_0(\refl)\subset T$ (respectively, $\Td^+_0(\refl)\subset \Td$) be
  as in Section~\ref{CCompactLie} (respectively, as above). It is easy to see that
  $\Td^+_0(\refl)$ is the group of $2$--primary torsion elements in
  $T^+_0(\refl)$, so that $T^+_0(\refl)$ is the topological closure of
  $\Td^+_0(\refl)$ in $T$, and $\Td^+_0(\refl)=\Td\cap T^+_0(\refl)$.
  Let $N(\refl)$ (respectively, $\Nd(\refl)$) be the centralizer of $T^+_0(\refl)$ in
  $N$ (respectively, the centralizer of $\Td^+_0(\refl)$ in $\Nd$). It follows
  that $\Nd(\refl)=\Nd\cap N(\refl)$. Let $\RG\refl$
  (respectively, $\RX\refl$) be the centralizer of $T^+_0(\refl)$ in $G$
  (respectively, the centralizer of $\Td^+_0(\refl)$ in $X$). There is a
  commutative diagram:
  \[
  \begin{CD}
     B \Td^-_0(\refl) @>>> BT^-_0(\refl) @. {} @. {} \\
      @VVV                @VVV            @.    @.  \\
     B\Nd(\refl) @>>>  BN(\refl) @>>> B\RG\refl @>>> B\RX\refl
  \end{CD}
  \]
  Again, $\Td^-_0(\refl)$ is the $2$--primary torsion subgroup of
  $T^-(\refl)$. Note that $\Nd(\refl)$ is a union of finite $2$--groups
  (Aside~\ref{HomotopicalConjugacy}). According to Proposition~\ref{SameCentralizers}
  and Proposition~\ref{SameConjugacies}, an element $x\in\Nd(\refl)$ is
  homotopically conjugate (Aside~\ref{HomotopicalConjugacy}) in $\RX\refl$
  to an element of $\Td^-_0(\refl)$ if and only if $x$ is
  algebraically conjugate in $\RG\refl$ to an element of
  $\Td^-_0(\refl)$. Since $x$ is of order a power of~$2$, this last
  occurs if and only if $x$ is conjugate in $\RG\refl$ to an element
  of $T^-_0(\refl)$. The lemma follows from the definitions of the
  root subgroups $N_\refl$ (Definition~\ref{DefineLGRootSubgroup}) and
  $\Nd_\refl$ (Definition~\ref{tDefineRootSubgroup}).
\end{proof}

\begin{proof}[Proof of Proposition~\ref{IdentifyNT}, Coxeter case]
  Let $(\Lc,W)$ be the \mrtlattice/ associated to $X$
  (Proposition~\ref{tInterpretMarkedTorus}), and $(L,W)$ a marked reflection
  lattice such that $(\Lc,W)\iso(\Ztwo\Tensor L, W)$. Associated to
  $(L,W)$ is a marked reflection torus $(T,W,\{h'_\refl\})$. The group
  $\Td\iso\Lc\Tensor\dTorus(1)$ is naturally isomorphic to the
  $2$--primary torsion subgroup of $T\iso L\Tensor\Torus(1)$, and under
  this identification, $h_\refl'=h_\refl$ for each reflection
  $\refl\in W$. Consequently, as in Section~\ref{CPassingToTwo}, the
  normalizer extension $\Normd(W)=\Normd(\Td,W,\{h_\refl\})$ is a
  subgroup of the extension $\Norm(W)= \Norm(T,W,\{h_\refl\})$
  (Definition~\ref{DefineNormalizerExt}). Choose a set $\{\srefl_i\}$ of simple
  reflections in $W$ (see Remark~\ref{CoxeterGeneration}).  By
  Section~\ref{GeneratorsAndRelations} and Definition~\ref{tDefineNormalizerExt},
  $\Normd(W)$ is generated by $\Td$, together with symbols $q_i$, one
  for each simple reflection $\srefl_i$, subject to the relations
  \eqref{FinalNTRelations} (with $T$ replaced by $\Td$).  As in the
  proof of Theorem~\ref{IdentifyTorusNormalizer}, in order to show that
  $\Normd(W)$ is isomorphic to $\Nd$, it is enough to find elements
  $x_i\in\Nd$, one for each simple reflection, such that $x_i$
  projects to $\srefl_i\in W$ and such that the $x_i$ satisfy
  relations \eqref{XRelations} (with $T$ replaced by $\Td$). For each
  simple reflection $\srefl_i$, let $\Nd_{\srefl_i}$ be the
  corresponding root subgroup, and $x_i$ some element in
  $\Nd_{\srefl_i} \setminus T^-_0(\srefl_i)$. We claim the elements
  $x_i$ satisfy the necessary relations. The first relation in
  \eqref{XRelations} is clear (Definition~\ref{DefineLXMarkings}), and the second
  follows from the fact that $x_i$ by construction projects to
  $\srefl_i\in W$. The final relation involves pairs $\{x_i,x_j\}$ of
  the generating elements.  Choose such a pair, let $A$ be the rank
  $r-2$ \ttorus/ which is the maximal divisible subgroup of $\Td$ on
  which $\srefl_i$ and $\srefl_j$ act as the identity, and let $X_A$
  be the centralizer of $A$ in $X$. Then $X_A$ is connected
  \cite[9.8]{rDWcenter} and of semisimple
  rank~$2$. In the notation of Lemma~\ref{CutDownToCentralizer}, $x_i$
  (respectively,  $x_j$) is in the root subgroup $\Nd_{A,\srefl_i}$
(respectively,
  $\Nd_{A,\srefl_j}$) of $\Nd$ determined by $X_A$, and so in order to
  verify the final relation between $x_i$ and $x_j$ we can replace $X$
  by $X_A$, or in particular, assume that $X$ has semisimple rank~$2$.
  In this case $X$ is the $2$--completion of a connected compact Lie group
  (Proposition~\ref{SemisimpleRankTwo}; note that this is a crucial step in the
  proof, which is due to Kitchloo and Notbohm), and the result essentially follows from
  Lemma~\ref{SameRootSubgroups} and the fact that the relation is satisfied
  in the Lie group case (Theorem~\ref{MainTitsTheorem}). The only remaining
  issue is to show that $\{\srefl_i,\srefl_j\}$ is a set of simple reflections
  for the Weyl group $W_A$ of $X_A$ (note that $W_A$ is the set of
  elements in $W$ which pointwise fix $A$). 

  Let $V$ be the
  codimension~$2$ subspace of $\R\Tensor L$ which is pointwise fixed
  by $\srefl_i$ and $\srefl_j$; it is easy to see that $W_A$ is the
  subgroup of $W$ consisting of element which fix $V$ pointwise. The
  group $W_A$ is generated by the reflections it contains (see
  \cite[Theorem~1.2(c)]{rHump}, or note that $W_A$ is the Weyl group
  of a connected $2$--compact group). The desired result now follows
  from \cite[Theorem~1.10(a)]{rHump}, given that, in the language of
  \cite{rHump}, a reflection in $W$ fixes $V$ pointwise if and only if
  the associated root is orthogonal to $V$, that is, it lies in the
  $\R$--span of the roots corresponding to $\srefl_i$ and $\srefl_j$.
\end{proof}

\begin{proof}[Proof of Proposition~\ref{IdentifyNT}, $\DI(4)$ case]
  Suppose that $X$ is of type $\DI(4)$.  Let $x$ be a nontrivial
  element of order $2$ in $\Td$, and $A$ the subgroup of $\Td$
  generated by $x$. Note that the Weyl group $\WD$ of $X$ acts
  transitively on the set of such elements $x$ (Remark~\ref{DIFourStuff}), so
  that the conjugacy class of $A$ as a subgroup of $X$ does not depend
  on the choice of $x$. Let $X_A$ be the centralizer of $A$ in $X$,
  and $W_A$ the group of elements in $\WD$ which fix $x$. The
  $2$--compact group $X_A$ is connected; one quick way to check this is
  to observe that all of the reflections in $\WD$ are nontrivial
  mod~$2$ (Remark~\ref{DIFourStuff}), and so by \cite[Section~7]{rDWcenter} $X_A$
  is connected if and only if $W_A$ is generated by reflections. But
  $W_A$ is generated by reflections; this can be verified directly,
  or by observing that in the particular case $X=\DI(4)$ \cite{rDWnew},
  the centralizer of such an $A$ is clearly connected, since it is the
  $2$--completion of $\operatorname{Spin}(7)$.

  By Remark~\ref{DIFourStuff} and \cite[Section~7]{rDWcenter}, the group $W_A$,
  which is the Weyl group of $X_A$, is isomorphic to $\{\pm1\}\times
  P(1,2)$, where $P(1,2)$ is the subgroup of $\GL(3,\Ftwo)$ given by
  matrices which agree with the identity in the first column. For
  cardinality reasons $W_A$ cannot contain $\WD$ as a factor, and so
  by Proposition~\ref{ClassifyCompleteLattices} $X_A$ is of Coxeter type. Since
  the index of $W_A$ in $\WD$ is~$7$, which is odd, the fact that
  Proposition~\ref{IdentifyNT} is true for $X$ follows from
  Remark~\ref{UseCentralizerCompatibility} and the fact proved above that
  Proposition~\ref{IdentifyNT} is true for $X_A$.
\end{proof}

\section{A centralizer lemma}\label{CCentralizer}

In this section we verify that the normalizer extension construction
given in Definition~\ref{tDefineNormalizerExt} 
behaves well in a certain sense when
it comes to taking centralizers of subgroups of the torus $\Td$.
Our proof of Proposition~\ref{IdentifyNT} in the $\DI(4)$ case depends on this
behavior.

We put ourselves in the context of Lemma~\ref{CutDownToCentralizer}; in
particular, $X$ is a connected $2$--compact group with marked
reflection torus $(\Td,W,\{h_\refl\})$ and discrete torus normalizer
$\Nd$, $A$ is a subgroup of $\Td$, $X_A$ is the centralizer of $A$ in
$X$ and $W_A$ is the subgroup of $W$ consisting of elements which
leave $A$ pointwise fixed \cite[7.6]{rDWcenter}. Recall that we are
assuming that $X_A$ is connected. According to Lemma~\ref{CutDownToCentralizer},
$(\Td, W_A, \{h_\refl\}_{\refl\in W_A})$ is the marked reflection
torus associated to $X_A$.  Let $\Norm(W)$ (respectively, $\Norm(W_A)$) denote
the normalizer extension of $W$ (respectively, $W_A$) obtained from $(\Td, W,
\{h_\refl\})$ (respectively, $(\Td, W_A,\{h_\refl\}_{\refl\in W_A})$.

\begin{lem}\label{CentralizerCompatibility}
  In the above situation, there is a commutative diagram:
  \[
    \begin{CD}
      1 @>>> \Td @>>> \Norm(W_A) @>>> W_A @>>> 1  \\
      @.    @V\|VV      @VVV        @V\subset VV     @.  \\
      1 @>>> \Td  @>>> \Norm(W)  @>>> W    @>>> 1
    \end{CD}
  \]
  In particular, the extension $\Norm(W_A)$ of $W_A$ is the pullback
  over $W_A\subset W$ of the extension $\Norm(W)$ of $W$.
\end{lem}

\begin{rem}\label{UseCentralizerCompatibility}
  Let $\Nd$ (respectively, $\Nd_A$) denote the normalizer of $\Td$ in $X$
  (respectively, $X_A$). It is clear from \cite[Section~7]{rDWcenter} that the
  extension $\Nd_A$ of $W_A$ is the pullback over $W_A\subset W$ of
  the extension $\Nd$ of $W$. We can therefore conclude from
  Lemma~\ref{CentralizerCompatibility} that if $\Nd\iso\Norm(W)$, then
  $\Nd_A\iso\Norm(W_A)$. More interestingly, suppose that the index of
  $W_A$ in $W$ is odd, so that (by a transfer argument) the
  restriction map $\HH^2(W;\Td)\to\HH^2(W_A;\Td)$ is injective and
  extensions of $W$ by $\Td$ are detected on $W_A$. In this case we
  can conclude that if $\Nd_A\iso\Norm(W_A)$, then $\Nd\iso\Norm(W)$.
\end{rem}

\begin{NumberedSubSection}{A double coset formula}\label{DiscussDoubleCoset}
The proof of Lemma~\ref{CentralizerCompatibility} depends upon a double
coset formula for Shapiro companions. If $u\co H\to G$ is an inclusion of
finite groups and $k\in \HH^2(H;\Z)$, write $u_{\#}(k)$ for the element
of $\HH^2(G;\Z[G/H])$ that corresponds to $k$ under Shapiro's lemma.
If $v\co K\to G$ is another subgroup inclusion, the problem is
to compute $v^*u_{\#}(k)\in\HH^2(K;\Z[G/H])$. Let $x_\alpha$ run
through a set of $(K,H)$ double coset representatives in $G$, so that
$G$ can be written as a disjoint union
\[
    G = \coprod_{\alpha} Kx_\alpha H.
\]
For each index $\alpha$, let $K_\alpha=K \cap x_\alpha H
x_\alpha^{-1}$, and note that there is a natural isomorphism of
$K$--modules
\begin{equation}\label{SumDecomposition}
     \Z[G/H] \iso \oplus_\alpha \Z[K/K_\alpha].
\end{equation}
Write $v_\alpha$ for the map 
\[
   v_\alpha\co  K_\alpha\to x_\alpha H
   x_\alpha^{-1}\RightArrow{x_\alpha^{-1}(\text{--})x_\alpha} H
\]
and $u^\alpha$ for the inclusion $K_\alpha\to K$. Then for each index
$\alpha$, $u^\alpha_\#\,v_\alpha^*(k)\in\HH^2(K;\Z[K/K_\alpha])$, so
that, under isomorphism \eqref{SumDecomposition}, the element $\oplus_\alpha
\,\,u^\alpha_\#\,v_\alpha^*(k)$ lies in $\HH^2(K;\Z[G/H])$.

\begin{lem}\label{DoubleCosetFormula}
In the above situation, for any $k\in \HH^2(H;\Z)$,
\[
     v^*u_\#(k)= \oplus_\alpha \,\,u^\alpha_\# \,v_\alpha^*(k).
\]
\end{lem}

\begin{proof}[Sketch of proof]
Consider the homotopy fibre square
\[
\begin{CD}
   \coprod_\alpha BK_\alpha @>\coprod Bv_\alpha >> BH\\
    @V\coprod Bu^\alpha VV                         @V Bu VV\\
    BK    @>Bv>>                  BG  
  \end{CD}
\]
in which the vertical maps are covering spaces with fibre $G/H$. An
element $k\in \HH^2(H;\Z)$ corresponds to an extension of $H$ by $\Z$
and thus to a fibration $p\co BE\to BH$ with fibre $B\Z=S^1$. The Shapiro
companion $u_\#(k)\in \HH^2(BG;\Z[G/H])$ similarly corresponds to some
fibration $(Bu)_\#(BE)$ over $BG$. One checks that $(Bu)_\#(BE)$ is
the fibration whose fibre over $x\in BG$ is the product, taken over
$y\in (Bu)^{-1}(x)$, of $p^{-1}(y)$. In other words, $u_\#$ is a kind
of geometric multiplicative transfer. By inspection, pulling $BE\to
BH$ back over $\coprod Bv_\alpha$ and applying the multiplicative
transfer $(\coprod Bu^\alpha)_\#$ gives the same result as pulling
$(Bu)_\#(BE)$ back over $Bv$. The lemma is the algebraic expression of
this.
\end{proof}
\end{NumberedSubSection}

\begin{proof}[Proof of Lemma~\ref{CentralizerCompatibility}]
Let $\refls$ be the set of reflections in $W$ and $\refls^A$ the set
of reflections in $W_A$, so that $\Z[\refls]$ is isomorphic as a
$W_A$--module to the sum $\Z[\refls^A]+\Z[\refls']$, where $\refls'$ is
the complement $\refls\setminus\refls^A$.
Because of the way in which the normalizer extension $\Norm(W_A)$
depends on the reflection extension $\RE(W_A)$, it is enough to check
that the restriction of the reflection extension $\RE(W)$ to $W_A$
is the sum of $\RE(W_A)$ and a semidirect product extension of $W_A$ by
$\Z[\refls']$. As in Section~\ref{tDiscussReflectionExt},
write $\refls=\coprod \refls_i$ as a union of $W$--conjugacy classes.
Each conjugacy class $\refls_i$ can then be written as a union
\[
  \refls_i = {\textstyle \coprod_j \refls_{i,j}}
\]
of $W_A$--conjugacy classes.  For each $i$ let $\crefl_i$ denote a
representative of the $W$--conjugacy class $\refls_i$, and
$\crefl_{i,j}$ a representative of the $W_A$--conjugacy class
$\refls_{i,j}$; the centralizer of $\crefl_i$ in $W$ is $C_i$ and the
isotropy subgroup of $\crefl_{i,j}$ in $W_A$ is $C_{i,j}$. The group
$C_i$ has a canonical extension by $\Z$ (\eqref{tCanonicalExtension});
if $\crefl_{i,j}\in W_A$ then $C_{i,j}$ is the centralizer of
$\crefl_{i,j}$ in $W_A$ and it too has a canonical extension by~$\Z$.
Denote the corresponding extension classes by $k_i\in\HH^2(C_i;\Z)$
and $k_{i,j}\in\HH^2(C_{i,j};\Z)$. Let $u^i\co C_i\to W$,
$u^{i,j}\co C_{i,j}\to W_A$, and $v\co W_A\to W$ be the inclusion maps. In
the notation of Section~\ref{DiscussDoubleCoset}, the task comes down to
checking for each index $i$ the validity of the formula
\begin{equation}\label{SoughtAfterFormula}
      v^* \,u^i_\#(k_i) = \oplus_{\crefl_{i,j}\in W_A} u^{i,j}_\#
      (k_{i,j}) \oplus_{\crefl_{i,j}\notin W_A} 0_{i,j},
\end{equation}
where $0_{i,j}$ denotes the zero element of $\HH^2(W_A;\Z[\refls_{i,j}])$.

The left hand side of \eqref{SoughtAfterFormula} can be evaluated with
Lemma~\ref{DoubleCosetFormula}.  The double cosets $W_A\backslash W/C_i$ in
question correspond to the orbits of the conjugation action of $W_A$
on $\refls_i$. For each such orbit $\refls_{i,j}$ with orbit
representative $\crefl=\crefl_{i,j}$, choose an element $t\in W$ such that
$t\crefl_i t^{-1}=\crefl$; the element $t$ is
then the double coset representative. Note that $C_{i,j}=W_A\cap
t C_it^{-1}$.  We have to consider the
maps
\[
\begin{aligned}
     v_{i,j}\co  C_{i,j} &\to tC_it^{-1}
      \RightArrow {t^{-1}(\text{--})t} C_i\\
      u^{i,j}\co C_{i,j}&\to W_A
\end{aligned}
\]
and compute $u^{i,j}_\#\,v_{i,j}^*(k_i)$. The key observation is that
if $\crefl$ is not contained in $W_A$, then $C_{i,j}$, which is
generated by the reflections it contains
\cite[1.5]{rSteinbergDifferential}, must lie in $tC_i\pperp t^{-1}$;
hence $v_{i,j}^*(k_i)=0$. If $\crefl$ is contained in $W_A$, then
$C_{i,j}$ is the centralizer of $\crefl$ in $W_A$, $ t^{-1}\crefl
t=\crefl_i$, and $t^{-1}C_{i,j}\pperp t\subset C_i\pperp$, so that
$v_{i,j}^*(k_i)=k_{i,j}$.
\end{proof}

\end{document}

%% file: gtoutput.tex
\def\ifplaintex{\expandafter\ifx\csname documentclass\endcsname\relax}

\ifplaintex 
\else
\fi

\expandafter\ifx\csname beginpicture\endcsname\relax
\expandafter\ifx\csname documentclass\endcsname\relax
\input pictex \else
\input prepictex \input pictex \input postpictex \fi\fi

\def\gt{{\mathsurround=0pt\it $\cal G\mskip-2mu$eometry \&\ 
$\cal T\!\!$opology}}        

\def\gtp{{\mathsurround=0pt\it $\cal G\mskip-2mu$eometry \&\ 
$\cal T\!\!$opology $\cal P\!$ublications}}  

\def\lognumber#1{\def\thelognumber{#1}}
\def\volumenumber#1{\def\thevolumenumber{#1}}
\def\papernumber#1{\def\thepapernumber{#1}}
\def\volumeyear#1{\def\thevolumeyear{#1}}

\def\pagenumbers#1#2{\def\startpage{#1}\def\finishpage{#2}}
\def\published#1{\def\publishdate{#1}}
\def\proposed#1{\def\theproposer{#1}}
\def\seconded#1{\def\theseconders{#1}}
\def\received#1{\def\receiveddate{#1}}
\def\revised#1{\def\reviseddate{#1}}
\def\accepted#1{\def\accepteddate{#1}}

\def\coverauthors#1{\def\thecoverauthors{#1}}
\def\asciiauthors#1{\def\theasciiauthors{#1}}
\def\asciiaddress#1{\def\theasciiaddress{#1}}
\def\asciiemail#1{\def\theasciiemail{#1}}

\long\def\asciiabstract#1{\long\def\theasciiabstract{#1}}

\let\\\par\let\thelognumber\relax
\let\thevolumenumber\relax\let\thepapernumber\relax
\let\thevolumeyear\relax\let\thesamplenumber\relax\let\startpage\relax
\let\finishpage\relax\let\publishdate\relax\let\receiveddate\relax
\let\reviseddate\relax\let\accepteddate\relax\let\theasciititle\relax
\let\theasciiauthors\relax\let\theasciiaddress\relax
\let\theasciiabstract\relax
\let\theasciiemail\relax\let\theshortauthors\relax\let\theshorttitle\relax
\let\thecoverauthors\relax

\long\def\maketitlep{   

\count0=\startpage

\gt\hfill      
\beginpicture
\setcoordinatesystem units <0.33truein, 0.33truein> point at 2.2 0.9
\setplotsymbol ({$\cal G$})
\plotsymbolspacing=9truept
\circulararc 315 degrees from 0 1 center at 0 0
\setplotsymbol ({$\cal T$})
\circulararc 315 degrees from 1 -1 center at 1 0
\endpicture
\break
{\small\ifx\thesamplenumber\relax 
Volume \else Sample
\fi\thevolumenumber\ (\thevolumeyear)
\startpage--\finishpage\nl
Published: \publishdate}
\vglue 0.5truein plus 0.4fil minus 0.1truein

{\parskip=0pt\leftskip 0pt plus 1fil\def\\{\par\smallskip}{\ifplaintex\large
\else\Large\fi\bf\thetitle}\par\medskip}   

\vglue 0pt plus 0.1fil 

{\parskip=0pt\leftskip 0pt plus 1fil\def\\{\par}{\sc\theauthors}
\par\medskip}

\vglue 0pt plus 0.1fil 

{\small\parskip=0pt\let\newline\\
{\leftskip 0pt plus 1fil\def\\{\par}{\sl\theaddress}\par}
\expandafter\ifx\theemail\relax    
\relax\else\vglue 5pt plus 0.02fil minus 2pt\def\\{\stdspace{\rm 
and}\stdspace} 
\cl{Email:\stdspace\tt\theemail}\fi
\ifx\theurl\relax                  
\relax\else\vglue 5pt plus 0.02fil minus 2pt\def\\{\stdspace{\rm 
and}\stdspace}
\cl{URL:\stdspace\tt\theurl}\fi\par}

\vglue 7pt plus 0.3fil minus 3pt

{\bf Abstract}
\vglue 5pt plus 0.1fil minus 2pt

\theabstract

\vglue 7pt plus 0.3fil minus 3pt

{\bf AMS Classification numbers}\quad Primary:\quad \theprimaryclass

Secondary:\quad \thesecondaryclass

\vglue 5pt plus 0.3fil minus 2pt

{\bf Keywords:}\quad \thekeywords

\vglue 10pt plus 0.5fil minus 5pt

{\small  Proposed: \theproposer\hfill Received: \receiveddate\nl
Seconded: \theseconders\hfill 
\ifx\reviseddate\relax                         
Accepted: \accepteddate                        
\else
Revised: \reviseddate                          
\fi}
\eject
}       

\let\maketitlepage\maketitlep
\let\maketitle\maketitlepage

\font\phead=cmsl9 scaled 950
\font=cmsl9 scaled 1050
\font\pnum=cmbx10 scaled 913
\font\lnum=cmbx10 
\font\pfoot=cmsl9 scaled 950
\font=cmsl9 scaled 1050
\ifplaintex
\headline{\vbox to 0pt{\vskip -4.5mm\line{\small\phead\ifnum
\count0=\startpage ISSN 1364-0380 (on line)
1465-3060 (printed) \hfill {\pnum\folio}\else\ifodd\count0\def\\{ }%
\ifx\theshorttitle\relax\thetitle\else\theshorttitle\fi\hfill{\pnum\folio}
\else\def\\{ and }{\pnum\folio}\hfill\ifx\theshortauthors\relax\theauthors
\else\theshortauthors\fi\fi\fi}\vss}}
\footline{\vbox to 0pt{\vglue 0mm\line{\small\pfoot\ifnum\count0=\startpage
\copyright\ \gtp\hfill\else
\gt, Volume \thevolumenumber\ (\thevolumeyear)\hfill\fi}\vss
}}
\else
\makeatletter
\def\@oddhead{{\small\lhead\ifnum\count0=\startpage ISSN 1364-0380 (on line)
1465-3060 (printed) \hfill {\lnum\number\count0}\else\ifodd\count0
\def\\{ }\ifx\theshorttitle\relax \thetitle \else\theshorttitle\fi\hfill
{\lnum\number\count0}\else\def\\{ and }{\lnum\number\count0}
\hfill\ifx\theshortauthors\relax 
\theauthors\else\theshortauthors\fi\fi\fi}}\def\@evenhead{\@oddhead}
\def\@oddfoot{\small\lfoot\ifnum\count0=\startpage\copyright\ \gtp\hfill\else
\gt, Volume \thevolumenumber\ (\thevolumeyear)\hfill\fi}
\def\@evenfoot{\@oddfoot}
\makeatother
\fi

\newwrite\gtoutfile
\long\gdef\makeheadfile{  
{\def\\{, }\def\s{ }
\immediate\openout\gtoutfile head.xxx
\immediate\write\gtoutfile{Proxy-for: \ifx\theasciiauthors\relax
\theauthors\else\theasciiauthors\fi\s<\ifx\theasciiemail\relax\theemail\else\theasciiemail\fi>}
\immediate\write\gtoutfile{\noexpand\\}
\immediate\write\gtoutfile{Authors: \ifx\theasciiauthors\relax
\theauthors\else\theasciiauthors\fi}
{\def\\{ }\immediate\write\gtoutfile{Title: \ifx\theasciititle\relax
\thetitle\else\theasciititle\fi}}
\immediate\write\gtoutfile{Subj-class: GT or SG or MG etc}
\immediate\write\gtoutfile{MSC-class: \theprimaryclass\ifx\thesecondaryclass\relax\else, \thesecondaryclass\fi}
\immediate\write\gtoutfile{Journal-ref: Geom. Topol. \thevolumenumber
(\thevolumeyear) \startpage-\finishpage}
\immediate\write\gtoutfile{Comments: Published by Geometry and Topology at}
\immediate\write\gtoutfile{\s\s http://www.maths.warwick.ac.uk/gt/GTVol\thevolumenumber/paper\thepapernumber.abs.html}
\immediate\write\gtoutfile{\noexpand\\}
\immediate\write\gtoutfile{}
\ifx\theasciiabstract\relax
\immediate\write\gtoutfile{\theabstract}\else
\immediate\write\gtoutfile{\theasciiabstract}\fi
\immediate\write\gtoutfile{}
\immediate\write\gtoutfile{\noexpand\\}
\immediate\write\gtoutfile{}
\immediate\closeout\gtoutfile}}  

\def\maketitlepage{\maketitlep\makeheadfile}
\let\maketitle\maketitlepage